\documentclass[review,sort&compress]{elsarticle}

\usepackage{lineno,hyperref}
\usepackage{amsfonts}
\usepackage{amsmath}
\usepackage{amssymb}
\usepackage{amsxtra}
\usepackage{stmaryrd}
\usepackage{mathabx}
\usepackage{pdflscape}
\usepackage{xcolor}
\usepackage{algorithm}
\usepackage{algorithmic}

\modulolinenumbers[5]

\journal{Journal of \LaTeX\ Templates}

    \newtheorem{theorem}{Theorem}[section]
    \newtheorem{lemma}[theorem]{Lemma}
    \newtheorem{proposition}[theorem]{Proposition}
    \newtheorem{corollary}[theorem]{Corollary}

    \newenvironment{proof}[1][Proof]{\begin{trivlist}
    \item[\hskip \labelsep {\bfseries #1}]}{\end{trivlist}}
    \newenvironment{definition}[1][Definition]{\begin{trivlist}
    \item[\hskip \labelsep {\bfseries #1}]}{\end{trivlist}}
    \newenvironment{example}[1][Example]{\begin{trivlist}
    \item[\hskip \labelsep {\bfseries #1}]}{\end{trivlist}}
    \newenvironment{notation}[1][Notation]{\begin{trivlist}
    \item[\hskip \labelsep {\bfseries #1}]}{\end{trivlist}}
    
    \newenvironment{remark}[1][Remark]{\begin{trivlist}
    \item[\hskip \labelsep {\bfseries #1}]}{\end{trivlist}}

\newcommand{\A}{\mathcal{A}}

\newcommand{\PP}{\mathcal{P}}
\newcommand{\PAB}{\PP_{[\alpha,\beta]}}
\newcommand{\CS}{\mathcal{S}}
\newcommand{\CI}{\mathcal{I}}
\newcommand{\pred}[1]{pred(#1)}

\newcommand{\myendproof}{$\boxempty$
\endproof}

\begin{document}

\begin{frontmatter}

\title{Intervals of Antichains and Their Decompositions}
%\tnoteref{mytitlenote}}

\author{Patrick De Causmaecker, Stefan De Wannemacker, Jay Yellen}
\address{Department of Computer Science, KU Leuven, KULAK, CODeS \& iMinds-ITEC}
\address{Department of Mathematics and Computer Science, Rollins College, Florida USA}


%\ead{patrick.decausmaecker@kuleuven-kulak.be}


\begin{abstract}
An antichain of subsets is a set of subsets such that no subset in the antichain is a proper subset of any other subset in the antichain.
The {\it Dedekind number} counts the total number of antichains of subsets of an $n$-element set.
This paper investigates the interval structure of the lattice of antichains.
Several partitioning theorems and counting formulas for the size of intervals are derived. 
%Properties of intervals in the lattice of antichains of subsets of a  universe of finite size are investigated.
%New objects and quantities in this lattice are defined.
%Expressions and numerical values are deduced for the number of connected antichains and the number of fully distinguishing antichains.
%The latter establish a connection with Stirling numbers of the second kind.
%Decomposition properties of intervals in the lattice of antichains are proven.
%A new operator allowing partitioning the full lattice in intervals derived from lower dimensional sub-lattices is introduced.
%Special posets underlying an interval of antichains are defined.
%The poset allows the derivation of a powerful formula for the size of an interval.
%This formula allows computing intervals in the six dimensional space.
%Combinatorial coefficients allowing another decomposition of the full lattice are defined.
%In some specific cases, related to connected components in graphs, these coefficients can be efficiently computed.
%This formula allows computing the size of the lattice of order 8 efficiently.
%This size is the number of Dedekind of order 8, the largest one known so far.
\end{abstract}

\begin{keyword}
%\texttt{elsarticle.cls}\sep \LaTeX\sep Elsevier \sep template
antichains, distributive lattice, posets, intervals, Dedekind numbers
%\MSC[2010] 00-01\sep  99-00
\end{keyword}

\end{frontmatter}

% \linenumbers

\section{Introduction}
\label{sec:introduction}
The {\it Dedekind number}, $M(n)$, counts the total number of antichains of subsets of an $n$-element set, or
equivalently, the number of monotonic Boolean functions on $n$ variables \cite{DEDEKIND,SLOANE,KLEITMAN_MARKOWSKY}. 
An antichain of subsets is a set of subsets such that no subset in the antichain is a proper subset of any other subset in the antichain. 
This collection of antichains forms a distributive lattice with the partial order for antichains, along with the {\it join} and {\it meet} operators, all defined in Section \ref{sec:preliminaries}.

In 1969, Kleitman \cite{KLEITMAN} obtained an upper bound for $log(M(n))$, which was improved by Kleitman and Markowsky \cite{KLEITMAN_MARKOWSKY} in 1975.
%\[ log |M(n)| < (1+O((log n)/n)) \binom{n}{\lfloor n/2 \rfloor}.\]
In 1981, Korshunov \cite{KORSHUNOV} used a sophisticated approach to give asymptotics for the Dedekind number itself. 
In 2001, Kahn \cite{KAHN} gave a simpler proof of the bound by Kleitman and Markowsky in  \cite{KLEITMAN_MARKOWSKY}. 

Finding a closed-form expression for $M(n)$, known as {\it Dedekind''s Problem}, has proven to be extremely difficult.  
Exact values have been found only for $n \leq 8$ \cite{CHURCH1940,CHURCH1965,WARD, WIEDEMANN}. 
This is sequence A000372 in Sloane's Online Encyclopedia of Integer Sequences \cite{SLOANE}. 
Recent attempts to improve the computational time are described in \cite{BAKOEV, FYDITEK}.
A related problem is counting the number of inequivalent antichains (i.e. antichains that do not transform in one another by a permutation of the elements) and a recent result was obtained in \cite{TAMON2014}. 

In this paper we investigate the interval structure of the lattice of antichains.
Intervals of antichains are defined according to their partial order.
We derive several partitioning theorems and counting formulas to compute the size of general intervals.

In addition to the basic definitions, notations, and conventions, Section \ref{sec:preliminaries} includes an efficient characterization of the join and meet operators, 
and introduces a new operator, the {\it direct product} of two antichains.  

Intervals of antichains are introduced in Section \ref{sec:intervalsofantichains},
and a key element of our counting strategy  - the {\it underlying poset} of an interval - is introduced. 
We give a complete characterization of such an underlying poset,
and illustrate its relationship to the interval structure.

Section \ref{sec:intervals} presents various decomposition properties of sets and intervals of antichains.  
We show when an interval can be expressed as a {\it direct join} of two other intervals,
and
how an antichain on a given set can lead to a partition of the interval of all antichains on that set.
The set of all antichains on a finite set can be partitioned into intervals using the direct product introduced in Section \ref{sec:preliminaries}. 

In Section \ref{sec:decomposebyproperty}, the size of the sets in antichains are used to find new decompositions. 
The concept of a {\it uniform} antichain is introduced referring to antichains where all of its subsets have the same size.
We derive a unique, nontrivial, decomposition of an antichain into uniform antichains.
This decomposition is the basis for the partition theorems that lead to
compact formulas for the number of antichains in an interval.

\section{Preliminaries}
\label{sec:preliminaries}
\begin{definition}
A set of subsets of a finite set is called an {\bf antichain} if no two distinct subsets are  {\it comparable}, i.e., neither is a subset of the other.  We will use lower-case Greek letters to denote the antichains. 
\end{definition}
\begin{notation}
\begin{itemize}
\item{
If $S \subset N= \{1,2,\ldots,n\}$, then the complete set of antichains in $S$, denoted by $\A_S$, is given by
\[
\A_S = \{ \alpha \subseteq 2^S |  \forall X,Y \in \alpha,   X \not= Y \Rightarrow X \not\subseteq Y \ and\  Y \not\subseteq X\}
\] 
We will use $\A_n \equiv \A_N$ interchangeably. 
% For $n=0$ we use $\A_0 = \{\{\},\{\emptyset\}\}$.
}
\item We denote the empty antichain $\{\}$ by  $\bot$, and the largest antichain $\{N\}$ by $\top$.
\end{itemize}
\end{notation}

\begin{definition}
$\A_S$ is a partially ordered set ({\em poset}), where the partial order $\le$ for any 
$\alpha, \beta \in \A_S$ is given by
$$\alpha \le \beta \Leftrightarrow \forall A \in \alpha, \exists B \in \beta, A \subseteq B$$  
\end{definition}

\begin{remark}
Throughout this paper, set relations and operators are used in different contexts.
For example, in the definition of the partial order $\le$, given above, set inclusion is applied to the individual elements of antichains, i.e., to subsets of $S$.
Elsewhere, set inclusion is applied to the antichains themselves.
In particular, $\alpha \subseteq \beta$ means that each subset in $\alpha$ is a subset in the antichain $\beta$.
Thus, $\alpha \subseteq \beta \Rightarrow \alpha \le \beta$, but the converse is false.
For example, if $\alpha = \{\{1\}\}$ and $\beta = \{\{1,2\}\}$, then $\alpha \le \beta$ but $\alpha \not\subseteq \beta$.
\end{remark}

\begin{definition}
Any two posets $(P_1,\le_1)$ and $(P_2,\le_2)$ are said to be {\bf isomorphic}, denoted $(P_1, \le_1) \cong (P_2, \le_2)$, if 
there exists an {\em order-preserving} bijection $f, P_1 \rightarrow P_2$, i.e.,
$\forall x,y \in P_1 : x \le_1 y \Leftrightarrow f(x) \le_2 f(y)$
\end{definition}

\begin{proposition}
If $S_1, S_2 \subseteq \{1,2,\ldots,n\}$ are of equal size, then $(\A_{S_1},\le) \cong  (\A_{S_2},\le)$, where $\le$ is the partial order for antichains.
\end{proposition}
\proof
Let $f:S_1 \rightarrow S_2$ be any bijection.
Clearly, for all $A,B \subseteq S_1$, $A \subseteq B \Leftrightarrow f(A) \subseteq f(B)$.
It follows that the mapping $\hat{f}$, given by $\hat{f}(\{X_1,X_2,\ldots,X_t\}) = \{f(X_1),f(X_2),\ldots,f(X_n)\}$ maps antichains in $\A_{S_1}$ to antichains in $\A_{S_2}$.
It remains to show that $\hat{f}$ preserves the antichain partial order $\le$.
That is, we must show $\forall \alpha,\beta \in \A_{S_1}$, $\alpha \le \beta \Leftrightarrow \hat{f}(\alpha) \le \hat{f}(\beta)$.
Suppose $\alpha \le \beta$ for $\alpha,\beta \in \A_{S_1}$, and let $Y \in \hat{f}(\alpha)$.
By definition of $\hat{f}$, $Y = f(X)$ for some $X \in \alpha$.
Since $\alpha \le \beta$, $X \subseteq W$ for some $W \in \beta$, and thus, $Y = f(X) \subseteq f(W) \in \hat{f}(\beta)$, which shows $\hat{f}(\alpha) \le \hat{f}(\beta)$.
A similar argument establishes the converse, which completes the proof. 
\myendproof

\begin{definition} Let $\alpha$ and $\beta$ be two antichains in $\A_N$. 
\begin{itemize} 
\item The {\bf join} of $\alpha$ and $\beta$, denoted $\alpha \vee \beta$, is the smallest antichain $\gamma$ such that $\alpha \le \gamma$ and $\beta \le \gamma$.
\item The {\bf meet} of $\alpha$ and $\beta$, denoted $\alpha \wedge \beta$, is the largest antichain $\gamma$ such that $\gamma  \le \alpha$ and $\gamma  \le \beta$.  
\end{itemize}
\end{definition}
\begin{definition} Let $\CS \subseteq 2^N$ be any collection of subsets of $N$. The {\bf maximum antichain} in $\CS$, denoted $maxAC(\CS)$, is the antichain that results from removing from $\CS$ all sets that are proper
subsets of some set in $\CS$. That is,
\[
maxAC(\CS) = \CS - \{A \in \CS | A \subsetneq B, \text{for some } B \in \CS\}
\]
\end{definition}

\begin{proposition} For any two antichains $\alpha$ and $\beta$,
\label{prop:max}
\begin{itemize}
\item[(i)] $\alpha \vee \beta =  maxAC(\alpha \cup \beta)$
\item[(ii)] $\alpha \wedge \beta =  maxAC(\{A \cap B | A \in \alpha, B \in \beta\})$
\end{itemize}
\end{proposition}
\proof
The two sets of inequalities, $maxAC(\alpha \cup \beta) \ge \alpha$ and  $maxAC(\alpha \cup \beta) \ge \beta$, follow directly from the definitions.
Thus, it remains to show that $maxAC(\alpha \cup \beta)$ is the smallest antichain satisfying its two inequalities, and $maxAC(\{A \cap B | A \in \alpha\text{ and } B \in \beta\})$
is the largest antichain satisfying its two inequalities.
\begin{itemize}
\item[(i)] Suppose $\gamma$ is any antichain such that $\gamma \ge \alpha$ and $\gamma \ge \beta$, and let $X \in maxAC(\alpha \cup \beta)$.
Then $X \in \alpha$ or $X \in \beta$, and in either case, $X \subseteq Y$ for some $Y \in \gamma$, which shows $maxAC(\alpha \cup \beta) \le \gamma$.
\item[(ii)] Now suppose $\gamma$ is any antichain such that $\gamma \le \alpha$ and $\gamma \le \beta$, and let $X \in \gamma$. 
Then $X \subseteq Y$ for some $Y \in \alpha$ and $X \subseteq Z$ for some $Z \in \beta$, and hence, $X \subseteq Y \cap Z$. By definition of the max operator, we have
$Y \cap Z \subseteq W$ for some $W \in maxAC(\{A \cap B | A \in \alpha\text{ and } B \in \beta\})$, which shows that $\gamma \le maxAC(\{A \cap B | A \in \alpha\text{ and } B \in \beta\})$.
\end{itemize}
This completes the proof.
\myendproof
Observe that for any two antichains,
$\alpha$ and $\beta$, $\alpha \cap \beta$ is an antichain, and $\alpha \cap \beta \le \alpha \wedge \beta$, whereas, $\alpha \cup \beta$ is not necessarily an antichain.

\begin{notation}
For any antichain $\alpha$, we use $\cup{\alpha}$ to denote $\cup_{X \in \alpha} X$.
\end{notation}

The following proposition is an immediate consequence of the definitions.  Parts (i) and (ii) establish that 
$(\A_N,\wedge,\vee)$ is a distributive lattice. 

\vskip .2cm
\begin{proposition}
\label{prop:distributive}
For any antichains $\alpha, \beta,  \gamma \in \A_N$, 

({\it i}) $\alpha \wedge (\beta \vee \gamma) = (\alpha \wedge \beta) \vee (\alpha \wedge \gamma)$

({\it ii}) $\alpha \vee (\beta \wedge \gamma) = (\alpha \vee \beta) \wedge (\alpha \vee \gamma)$

({\it iii}) $\alpha \le \beta \implies (\alpha \wedge \gamma) \le (\beta \wedge \gamma)$
\end{proposition}

We end this section by introducing a new operator which will show very powerful. It will be used later on, a.o., in Proposition \ref{prop:intervalredundancy} 
and most importantly for decomposition in Theorem   \ref{the:intervalintervaldecomposition}.
\begin{definition}
Let $\alpha$ and $\beta$ be two antichains satisfying $(\cup{\alpha} )\cap (\cup{\beta}) = \emptyset$.
The {\bf direct product} $\alpha \otimes \beta$ is the antichain given by
\[
\alpha \otimes \beta = \{A \cup B | A \in \alpha, B \in \beta\}
\]
\end{definition}

\section{Intervals of Antichains}
\label{sec:intervalsofantichains}

\begin{definition}
For $\alpha, \beta \in \A_N$, the (closed) {\bf interval} from $\alpha$ to $\beta$ is defined by
$[\alpha,\beta] = \{\chi \in \A_N| \alpha \le \chi \le \beta\}$.  $\alpha$ and $\beta$ are called the {\bf bottom} and {\bf top} of the interval.  The half-open and open intervals, 
$(\alpha,\beta]$, $[\alpha,\beta)$, and $(\alpha,\beta)$ use strict inequalities, accordingly.
\end{definition}

Observe that $[\alpha,\beta] \not= \emptyset$ if and only if $\alpha \le \beta$.
Moreover, if $\chi_1, \chi_2 \in [\alpha,\beta]$, then $\chi_1 \wedge \chi_2 \in [\alpha,\beta] $ and $\chi_1 \vee \chi_2 \in [\alpha,\beta]$. 

The next two propositions are also immediate consequences of the definitions.

\begin{proposition}
\label{prop:intervalintersect}
For any four antichains $\alpha_1, \alpha_2, \beta_1, \beta_2$,
$$[\alpha_1,\beta_1] \cap [\alpha_2,\beta_2] = [\alpha_1 \vee \alpha_2,\beta_1 \wedge \beta_2]. $$
\end{proposition}

\begin{proposition}
For $\alpha \le \beta$, the mapping $h_{[\alpha,\beta]} : \A_N \rightarrow [\alpha,\beta]$ defined by  $h_{[\alpha,\beta]}(\chi) = \alpha  \vee (\chi \wedge \beta )$ is a lattice homomorphism.  That is, for any $\chi_1, \chi_2 \in \A_N$,
\begin{eqnarray*}
h_{[\alpha,\beta]}(\chi_1\vee \chi_2)=h_{[\alpha,\beta]}(\chi_1)\vee h_{[\alpha,\beta]}(\chi_2)\\
h_{[\alpha,\beta]}(\chi_1\wedge \chi_2)=h_{[\alpha,\beta]}(\chi_1)\wedge h_{[\alpha,\beta]}(\chi_2)\\
\end{eqnarray*}
\end{proposition}

\begin{definition}
Let $\alpha, \beta \in \A_{N}$ be two antichains with $\alpha \le \beta$.  
\begin{itemize}
\item{ The {\bf underlying poset of the interval $I = [\alpha,\beta]$}, denoted $\PP_I$ or $\PAB$, is the poset of subsets of $N$ under set inclusion given by 
$\PAB = \{X \subseteq N | \alpha \vee \{X\} \in (\alpha,\beta]\}$ 
or, equivalently, 
$\PAB = \{X \subseteq N | \{X\} \not\le \alpha \text{ and } \{X\} \le \beta]\}.$
}
\item {
For a set of sets, $\CS \subseteq 2^N$, the poset $(\CS,\subseteq)$ {\bf spans} the interval $[\alpha,\beta]$ if for each antichain $\gamma \in [\alpha,\beta]$, there exists an antichain $\gamma' \subseteq \CS$ such that 
$\gamma = \alpha \vee \gamma'$.
}
\item{
A poset $(\CS, \subseteq)$ is called an {\bf interval poset} if there exists an interval 
$[\alpha, \beta]$ for some $\alpha, \beta \in \A_N$ such that $\CS = \PAB$ and the interval spanned by an interval poset $(\CS,\subseteq)$, is denoted by $\CI_{\CS}$,
and hence
$
\CS = \PP_{\CI_{\CS}}
$
}
\end{itemize}
\end{definition}
\begin{remark}
When referring to the poset $\CS$ under set inclusion, we will not distinguish between $\CS$ and $(\CS,\subseteq)$.
\end{remark}

\begin{proposition}
\label{prop:posetspan}
$\PAB$ spans the interval $[\alpha,\beta]$.
\end{proposition}
\proof
Let $\gamma \in [\alpha,\beta]$.  
The result is trivially true if $\gamma = \alpha$, by letting $\gamma' = \bot$. Assume that $\gamma > \alpha$, and let $\gamma' = \gamma - \alpha$.  Since $\gamma'$ is an antichain and $\gamma = \alpha \vee \gamma'$, it suffices to show that $\gamma' \subseteq \PAB$.  Let $X \in \gamma'$, and consider the antichain $\alpha \vee \{X\}$.  If $\alpha \vee \{X\} = \alpha$, then 
$\{X\} \le \alpha$.  But $X \notin \alpha \implies \exists A \in \alpha$ such that $X \subsetneq A$.  Moreover, $\alpha \le \gamma \implies \exists C \in \gamma$ such that $A \subseteq C \implies X \subsetneq A \subseteq C$, which would contradict the antichain property of $\gamma$ ($X,C \in \gamma$, with $X \neq C$, $X \neq \emptyset$, but $X \subset C$).  Thus, $\alpha \vee \{X\}  > \alpha$, and hence, $X \in \PAB$, which completes the proof.
\myendproof

\begin{corollary}
The poset of {\it all} subsets of $N$, $(2^N,\subseteq)$, spans $\A_N$.
\end{corollary}
\proof 
$(2^N,\subseteq) = \PP_{[\bot,\top]}$ spans $[\bot,\top] = \A_N$. 
\myendproof

\begin{definition} 
For $X \subseteq N$, the antichain of {\it immediate subsets} of X, denoted $pred(X)$, is given by  
\[
pred(X) = \{X - \{x\} | x \in X\}
\]
%\item For an antichain $\alpha$, we let $\alpha^- = \bigvee_{X \in \alpha} \pred{X}$ and $\alpha^+ = \bigvee_{X \in 2^N, \pred{X} \le \alpha} \{X\}$.\break  
%Observe that $\emptyset^- =  \{\emptyset\}^- = \bot$, and $\bot^+ = \{\emptyset\}$.
\end{definition}

The following theorem characterizes interval posets and shows how to compute the bottom and top of their spanned intervals.

\begin{theorem}
\label{the:posetcharacteristic}
A collection of subsets $\CS \subseteq 2^N$ is an interval poset if and only if
\[
[*]\text{              } \forall A_1,A_2 \in \CS, \forall C  \subseteq N , A_1 \subseteq C  \subseteq A_2 \Rightarrow C \in \CS
\]
Moreover, if $\CS$ is an interval poset, then the interval it spans is
\[
\CI_{\CS} = [\bigvee_{X \in \CS} ((\pred{X}) - \CS),\bigvee_{X  \in \CS} \{X\}]
\]
\end{theorem}
\proof
($\Longrightarrow$)
Suppose $\CI_{\CS} = [\bigvee_{X \in \CS} ((\pred{X}) - \CS),\bigvee_{X  \in \CS} \{X\}]$, or equivalently, $\CS = \PP_{[\bigvee_{X \in \CS} ((\pred{X}) - \CS),\bigvee_{X  \in \CS} \{X\}]}$.
To show that the condition [*] holds, let $A_1 \subseteq C \subseteq A_2$ where $A_1,A_2 \in \CS$ and $C \subseteq N$.
We must show $C \in \CS$.  Since $A_1 \in \CS$, $\{A_1\}\not\le \bigvee_{X \in \CS} ((\pred{X}) - \CS)$.  But $A_1 \subseteq C \implies \{A_1\} \le \{C\}$, and hence, 
$\{C\}\not\le \bigvee_{X \in \CS} ((\pred{X}) - \CS)$.  Since $A_2 \in \CS$ and $C \subseteq A_2$, we have 
$\{C\} \le \{A_2\} \le \bigvee_{X  \in \CS} \{X\}$, and hence, $C \in \CS$, which shows that the condition [*] holds.

\noindent ($\Longleftarrow$)  Suppose that the condition [*] holds. We must show 
\[
\CS = \PP_{[\bigvee_{X \in \CS} ((\pred{X}) - \CS),\bigvee_{X  \in \CS} \{X\}]}
\]

\noindent ($\subseteq$) Let $A \in \CS$.  We have $\{A\} \le \bigvee_{X  \in \CS} \{X\}$, and it remains to show 
$\{A\}\not\le \bigvee_{X \in \CS} ((\pred{X}) - \CS)$.
If $\{A\}\le \bigvee_{X \in \CS} ((\pred{X}) - \CS)$, then  
$A \subseteq B$ for some $B \in \bigvee_{X \in \CS} ((\pred{X}) - \CS)$.  Thus, $B = X - \{x\} \not\in \CS$ for some $X \in \CS$.  But the condition implies that $B \in \CS$ since $A \subseteq B \subseteq X$.  This contradiction completes the proof that $\CS \subseteq \PP_{[\bigvee_{X \in \CS} ((\pred{X}) - \CS),\bigvee_{X  \in \CS} \{X\}]}$.

\noindent ($\supseteq$) Let $C \in \PP_{[\bigvee_{X \in \CS} ((\pred{X}) - \CS),\bigvee_{X  \in \CS} \{X\}]}$.  
We must show $C \in \CS$.\hfill\break  
Since $\{C\} \le \bigvee_{X  \in \CS} \{X\}$, $C \subseteq D$ for some $D \in \bigvee_{X \in \CS} \{X\}$, and by definition of the join, $D \in \CS$. Thus, $C \subseteq X$ for at least one $X \in \CS$.
Among all such $X$, let $B$ be a minimal one.  If $C=B$, then $C \in \CS$.  Otherwise, $C \subseteq B - \{b\}$ for some 
$b \in B$.  By the minimality of $B$, $B - \{b\} \not\in \CS$ and $B - \{b\} \in (\pred{b} - \CS)$. Thus, $\{C\} \le \bigvee_{X \in \CS} ((\pred{X}) - \CS)$. By definition, no element $C \in \PP_{[\bigvee_{X \in \CS} ((\pred{X}) - \CS),\bigvee_{X  \in \CS} \{X\}]}$ can satisfy $C \le \bigvee_{X \in \CS} ((\pred{X}) - \CS)$, and hence, we have reached a contradiction. This completes the proof. 
\myendproof

%
% Poset Decomposition
%
Lemma \ref{lem:posetredundancy} shows how to remove specific elements from the universe without changing the structure of the interval poset.
%
% Poset Redundancy
%
\begin{lemma}
\label{lem:posetredundancy}
Suppose that $\CS$ is a poset for which there exists a set $A$ such that  $A \subseteq X$ for all $X \in \CS$. Then the poset $\{X - A| X \in \CS\}$ is isomorphic to $\CS$.
\end{lemma}
\proof
The mapping $f : \CS \rightarrow \{X - A| X \in \CS\}$, given by $f(X) = X - A$, is easily shown to be an order-preserving bijection.
\myendproof

\begin{proposition}
\label{prop:isomorphicintervals}
If $\CS_1$ and $\CS_2$ are isomorphic posets, then the intervals they span are isomorphic. That is, $\CI_{\CS_1} \cong \CI_{\CS_2}$.
\end{proposition}
\proof
By Theorem \ref{the:posetcharacteristic}, 
$
\CI_{\CS_i} = [\bigvee_{X \in \CS_i} ((\pred{X}) - \CS_i),\bigvee_{X  \in \CS_i} \{X\}], i = 1,2
$.
Let $f : \CS_1 \rightarrow \CS_2$ be an order-preserving bijection.
For each $\chi \in \CI_{\CS_1}$, we have $\chi = [\bigvee_{X \in \CS_1}((\pred{X}) - \CS_1)] \vee \chi_1$ for some antichain $\chi_1 \subseteq \CS_1$.
Consider the mapping $\hat{f} : \CI_{\CS_1} \rightarrow \CI_{\CS_2}$ given by
\[
\hat{f}(\chi) = [\bigvee_{X \in \CS_2}((\pred{X}) - \CS_2)] \vee f(\chi_1)
\]
$\hat{f}$ is easily shown to be order preserving.
To show $\hat{f}$ is a bijection, suppose $\hat{f}(\chi) = \hat{f}(\gamma)$ for $\chi, \gamma \in I_{\CS_1}$.
Then \[ [\bigvee_{X \in \CS_2}((\pred{X}) - \CS_2)] \vee f(\chi_1) =  [\bigvee_{X \in \CS_2}((\pred{X}) - \CS_2)] \vee f(\gamma_1)\] where $\chi_1, \gamma_1 \in \CS_1$ .
We must show that $\chi = \gamma$, i.e., \[ [\bigvee_{X \in \CS_1}((\pred{X}) - \CS_1)] \vee \chi_1)=  [\bigvee_{X \in \CS_1}((\pred{X}) - \CS_1)] \vee \gamma_1\]
Since $f$ is a bijection, it suffices to show that $f(\chi_1) = f(\gamma_1)$.
Let $A \in f(\chi_1)$.  We have $A \not\in [\bigvee_{X \in \CS_2}((\pred{X}) - \CS_2)] $ (since $f(\chi_1) \subseteq \CS_2$).
But $A \in  [\bigvee_{X \in \CS_2}((\pred{X}) - \CS_2)] \vee f(\gamma_1) \Rightarrow A \in f(\gamma_1)$, and hence, $f(\chi_1) \subseteq f(\gamma_1)$.
Similarly,  $f(\gamma_1) \subseteq f(\chi_1)$, which completes the proof.
\myendproof

%
%
% Interval Redundancy
%
Proposition \ref{prop:intervalredundancy} builds on Lemma \ref{lem:posetredundancy} to associate a large set of intervals with any interval for which the span of the top is a true subset of $N$:
\begin{proposition}
\label{prop:intervalredundancy}
Given two antichains $\chi' \le \chi$ in $\A_N$ such that $N - \cup_{\chi} \not= \emptyset$.
All intervals of the form $[(\{A\} \otimes \chi') \vee (\pred{A} \otimes \chi),\{A\} \otimes \chi]$ where $A \subseteq N - \cup_{\chi}$, are isomorphic to $[\chi',\chi]$. 
\end{proposition} 
\proof
Let $\alpha = (\{A\} \otimes \chi') \vee (\pred{A} \otimes \chi)$ and $\beta = \{A\} \otimes \chi$.
Any element  $X \in P_{[\alpha,\beta]}$ must be a subset of one of the elements of $\beta$, so let 
$X \subseteq A \cup B$ with $B \in \chi$. 
If $A \not\subseteq X$, then $X \subseteq A' \cup B$ for some $A' \in \pred{A}$ and 
$\alpha \vee \{X\} = \alpha$, so that $X \not\in P_{[\alpha,\beta]}$.
So, $P_{[\alpha,\beta]}$ satisfies the condition of Lemma \ref{lem:posetredundancy}, and it is isomorphic to the poset $\{X - A | X \in P_{[\alpha,\beta]}\} $.
We now find
\[
P_{[\alpha,\beta]} = \{A \cup C | \exists B \in \chi: C \subseteq B \ and\ \not\exists B' \in \chi': C \subseteq B'\} \text{ and}
\]
\[
\{X - A | X \in P_{[\alpha,\beta]}\} = \{C | \exists B \in \chi: C \subseteq B \ and\ \not\exists B' \in \chi': C \subseteq B'\} = P_{[\chi',\chi]}
\]
Since isomorphic posets span isomorphic intervals (Proposition \ref{prop:isomorphicintervals}), this completes the proof.
\myendproof

\section{Interval Decomposition}
\label{sec:intervals}
This section presents two techniques for the decomposition of intervals.
In a first subsection, we use unique decomposition of antichains as the join of two or more other antichains.
In the second subsection, we study partitions intervals as the union of disjoint intervals.
\subsection{Decomposition of an Interval as a Direct Join}
\begin{definition}
Let $[\alpha_1,\beta_1]$ and $[\alpha_2,\beta_2]$ be any two intervals of antichains. Their {\bf join}, denoted $[\alpha_1,\beta_1] \vee [\alpha_2,\beta_2]$, is given by
\[
[\alpha_1,\beta_1] \vee [\alpha_2,\beta_2] = \{\chi| \chi = \chi_1 \vee \chi_2 \text{ for some } \chi_i \in [\alpha_i,\beta_i] \}
\]
It is easy to show that $[\alpha_1,\beta_1] \vee [\alpha_2,\beta_2] = [\alpha_1 \vee \alpha_2, \beta_1 \vee \beta_2]$.
The join $[\alpha_1,\beta_1] \vee [\alpha_2,\beta_2]$ is called a {\bf direct join}, and is denoted as $[\alpha_1,\beta_1] \ovee  [\alpha_2,\beta_2]$,
if the decomposition for each $\chi = \chi_1 \vee \chi_2$ is always unique.
\end{definition}

Proposition \ref{prop:posetdecomposition} shows how interval posets can be decomposed if the respective constituents are incomparable.
\begin{proposition}
\label{prop:posetdecomposition}
Suppose that $\CS_1$ and $\CS_2$ are two interval posets such that no two sets $X_1 \in \CS_1$ and $X_2 \in \CS_2$ are in a subset/superset relationship.
Then the poset $\CS_1 \cup \CS_2$ is an interval poset, and we have $\CI_{\CS_1 \cup \CS_2} = \mathcal{I}_{\CS_1} \ovee \CI_{\CS_2}$. 
% \proof
% Immediate
% \myendproof
\end{proposition}
\proof
To show condition [*] of Theorem \ref{the:posetcharacteristic} holds, let $A, B \in \CS_1 \cup \CS_2$, and suppose $A \subseteq C \subseteq B$.
We must show $C \in \CS_1 \cup \CS_2$.
Since $A \subseteq B$, both $A$ and $B$ must be in the same $\CS_i$.
It follows (by condition [*]) that $C \in \CS_i$, and hence $C \in \CS_1 \cup \CS_2$.
\myendproof

The next theorem allows decompositions of intervals as direct joins. Its proof uses the following lemma.
\begin{lemma}
\label{lem:directjoindecomposition}
For any antichains $\beta$ and $\delta$, if $\{X\} \not\le \beta$, then
\[
X \in \delta \Leftrightarrow X \in \beta \vee \delta
\]
\end{lemma}
\proof
This follows directly from Proposition \ref{prop:max}$(i)$.
\myendproof
%
% Interval Decomposition
%
\begin{theorem}
\label{the:directjoindecomposition}
Let $\alpha, \nu_1, \nu_2$ be any antichains such that $\nu_1 \wedge \nu_2 \le \alpha \le \nu_1 \vee \nu_2$.
Then the interval $[\alpha,\nu_1 \vee \nu_2]$ has the following direct-join decompositions:
\[
[\alpha,\nu_1 \vee \nu_2] = [\alpha,\alpha \vee \nu_1] \varovee [\alpha,\alpha \vee \nu_2] = [\alpha \wedge \nu_1,\nu_1] \varovee [\alpha \wedge \nu_2,\nu_2]
\]
\end{theorem}
\proof
\begin{itemize}
\item[(i)] Show $[\alpha,\nu_1 \vee \nu_2] = [\alpha,\alpha \vee \nu_1] \vee [\alpha,\alpha \vee \nu_2]$.\\
Suppose $\chi \in [\alpha,\nu_1 \vee \nu_2]$, and let $\chi_i = (\chi \wedge \nu_i) \vee \alpha, i=1,2$. Then
\begin{eqnarray*}
\chi_1 \vee \chi_2 & = & [(\chi \wedge \nu_1) \vee \alpha] \vee [(\chi \wedge \nu_2) \vee \alpha] \\
& = & (\chi \wedge \nu_1) \vee (\chi \wedge \nu_2) \vee \alpha \\
& = & [\chi \wedge (\nu_1 \vee \nu_2)] \vee \alpha \\
& = & \chi \vee \alpha = \chi
\end{eqnarray*}
Moreover, $\chi_i = \alpha \vee (\chi \wedge \nu_i) \le \alpha \vee \nu_i$, and hence, $\chi_i \in [\alpha,\alpha \vee \nu_i], i= 1,2$,
which shows $[\alpha,\nu_1 \vee \nu_2] \subseteq [\alpha,\alpha \vee \nu_1] \vee [\alpha,\alpha \vee \nu_2]$.
For the reverse inclusion, we have
\[
\chi \in [\alpha,\alpha \vee \nu_i] \Rightarrow \alpha \le \chi \le \alpha \vee \nu_i \Rightarrow \alpha \le \chi \le \nu_1 \vee \nu_2 
\Rightarrow \chi \in [\alpha, \nu_1 \vee \nu_2]
\]
which completes the proof of (i).
\item[(ii)] Show that the join is a direct join by establishing uniqueness.\\
Suppose $\chi_1 \vee \chi_2 = \gamma_1 \vee \gamma_2$, where $\chi_i,\gamma_i \in [\alpha,\alpha \vee \nu_i], i=1,2$.
We must show $\chi_i = \gamma_i, i=1,2$.
\end{itemize}
Let $X \in \chi_1$. We must show  $X \in \gamma_1$.\\
{\bf Case 1:} $\{X\} \not\le \alpha$.\\
Then $\{X\} \not\le \chi_2$, since otherwise, 
\[
\{X\} \le \chi_1 \wedge \chi_2 \le (\alpha \vee \nu_1) \wedge (\alpha \vee \nu_2) 
= \alpha \vee (\nu_1 \wedge \nu_2)  = \alpha
\]
which would contradict the hypothesis of Case 1.
% Similarly, $\{X\} \not\le \gamma_2$.
Thus, by Lemma \ref{lem:directjoindecomposition}, $X \in \chi_1 \vee \chi_2 = \gamma_1 \vee \gamma_2$.
But $\chi_1 \wedge \gamma_2 \le \alpha$ (arguing as above), and hence, $\{X\} \not\le \gamma_2$.
It follows that $X \in \gamma_1$, by Lemma \ref{lem:directjoindecomposition}.\\
{\bf Case 2:} $\{X\} \le \alpha$.\\
Then $X \subseteq Y$ for some $Y \in \alpha$ and $Y \subseteq Z$ for some $Z \in \chi_1$ (since $\alpha \le \chi_1$).
Since $\chi_1$ is an antichain, it follows that $X = Y = Z$ and hence $X \in \alpha$.
Since $\alpha \le \gamma_1$, $X \subseteq C$ for some $C \in \gamma_1$.
If $X = C$ then $X \in \gamma_1$, which would complete Case 2. If $X \subsetneq C$, then % $\{X\} < \{C\} \Rightarrow
$ \{C\} \not\le \chi_1$ (since $\chi_1$ is an antichain). Moreover, since $\{C\} \le \gamma_1 \vee \gamma_2 = \chi_1 \vee \chi_2$, we have $\{C\} \le \chi_2$ (by Lemma \ref{lem:directjoindecomposition}).
Thus, $X \subsetneq C \subseteq D$ for some $D  \in \alpha$, which contradicts the antichain property.
This shows that $X = C$, and hence $X \in \gamma_1$, which completes the proof that $\chi_1 \subseteq \gamma_1$.\\

\noindent The three other subset inclusions, $\gamma_1 \subseteq \chi_1$, $\chi_2 \subseteq \gamma_2$, and $\gamma_2 \subseteq \chi_2$ follow by similar arguments.

\noindent A similar argument can be used to show
\[
[\alpha,\nu_1 \vee \nu_2] = [\alpha \wedge \nu_1,\nu_1] \varovee [\alpha \wedge \nu_2,\nu_2]
\]
which completes the proof of Theorem \ref{the:directjoindecomposition}.
\myendproof

\subsection{Partitioning an Interval into Disjoint Intervals}
\label{sec:partitioningbyintervals}
\label{sub:antichainpartition}

\begin{definition} For any antichain $\chi \in \A_N$, let $\widecheck\chi$, called the {\bf largest nondominating antichain of $\chi$}, denote the largest antichain with the property that none of its elements is a superset of any element in $\chi$.  
Observe that the set of all such antichains, $\{\alpha \in \A_N | \forall X \in \chi, \{X\}\not\le \alpha\}$, is closed under the join operator, and hence, has a largest element. 
\end{definition}

\begin{proposition} $\forall A \subseteq N, \hskip.1cm  \widecheck{\{A\}} = \{N - \{a\} | a \in A\}$.
\end{proposition}
\proof
Suppose $X \in \widecheck{\{A\}}$. By definition, $\hskip.1cm A \not\subseteq X$, which implies $X\subseteq N - \{a\}$ for some $a \in A$. 
Thus, $\widecheck{\{A\}} \le \{N - \{a\} | a \in A\}$.  Conversely, 
if $X \in \{N - \{a\} | a \in A\}$,  then $A \not\subseteq X$, which implies 
$\{N - \{a\} | a \in A\} \le \widecheck{\{A\}}$ (since $\widecheck{\{A\}}$ is the largest antichain with that property).
\myendproof

%In Theorem \ref{the:largestnondominatingdecomposition} (and Corollary \ref{cor:largestnondominatingdecomposition}) a partition is associated with each antichain in (an interval of)  $\A_N$.
%To ease notation, we introduce the symbol $\check\chi$ for the largest antichain that does not dominate any of the sets in an antichain $\chi \in \A_N$:
%\begin{eqnarray}
%\label{eq:largestnondominating}
%\check\chi & \in & \A_N \nonumber \\
%\forall X \in \chi & : & \{X\} \not\le \check{\chi} \nonumber \\
%\forall \rho \in \A_N : (\forall X \in \chi : \{X\} \not\le \rho) & \Rightarrow & \rho \le \check{\chi}
%\end{eqnarray}
%This maximum exists since the criterion is closed under the join operation and computing the largest non dominating antichain $\widecheck{(.)}$ is straightforward:
%\begin{eqnarray*}
%\forall A \subseteq N : \widecheck{\{A\}} & = & \{N \backslash \{a\} | a \in a\}\\
%\forall \alpha,\beta \in  \A_N: \widecheck{\alpha \vee \beta} & = & \widecheck{\alpha} \wedge \widecheck{\beta}
%\end{eqnarray*}
%We can now prove the following

\begin{lemma}
\label{lem:disjointness}
Let $\alpha, \chi \in \A_N$, where $\chi \subseteq \alpha$.  
For each antichain $\sigma \in [\chi,\widecheck{(\alpha - \chi)}]$, $\chi = (\sigma \wedge \alpha)\cap \alpha$.
\end{lemma}
\proof
Since $\chi \le \sigma$ and $\chi \subseteq \alpha$, we have
$\chi = \chi \wedge \alpha \le \sigma \wedge \alpha \le (\sigma \wedge \alpha)\cap \alpha$.
To establish the reverse inequality, suppose that $A \in (\sigma \wedge \alpha)\cap \alpha$.
Then $\{A\} \le \sigma$ and $A \in \alpha$, and since $\sigma \in [\chi,\widecheck{(\alpha - \chi)}]$, we have $\{A\} \le \widecheck{(\alpha - \chi)}$.  It follows that $A \in \chi$, since otherwise, $A \in \alpha - \chi$ would contradict the definition of $\widecheck{(\alpha - \chi)}$. 
\myendproof

\begin{theorem}
\label{the:largestnondominatingdecomposition}
For any antichain $\alpha \in \A_N$, 
\[
 \A_N = \bigcup_{\chi \subseteq \alpha} [\chi,\widecheck{(\alpha - \chi)}]
\]
where the intervals are pairwise disjoint.
\end{theorem}
\proof
The pairwise disjoint property follows from Lemma \ref{lem:disjointness}.
For the set equality, it suffices to show 
$\A_N \subseteq \bigcup_{\chi \subseteq \alpha} [\chi,\widecheck{(\alpha - \chi)}]$.
Suppose that $\sigma \in  \A_N$, and let $\chi = (\sigma \wedge \alpha) \cap \alpha$.
Clearly, $\chi \le \sigma$.  To show $\sigma \le \widecheck{(\alpha - \chi)}$, we first show that $\sigma$ does not dominate any set in $\alpha - \chi$, i.e., $A \in \alpha - \chi \implies \{A\} \not\le \sigma$.  Let 
$A \in \alpha - \chi$, and suppose that $\{A\} \le \sigma$.  Since $A \in \alpha$, we have 
$\{A\} \le \sigma \wedge \alpha$, which implies that $A \subseteq B$ for some 
$B \in \sigma \wedge \alpha$.  But $\sigma \wedge \alpha \le \alpha \implies B \subseteq C$ for some $C \in \alpha$, and hence, $A\subseteq B\subseteq C$, where $A, C\in\alpha$.  Since $\alpha$ is an antichain, we must have $A=B=C$.  Thus, $A\in \sigma\wedge\alpha$, which implies $A\in \chi$, a contradiction, showing that $\sigma$ does not dominate any set in $\alpha - \chi$.  It follows that 
$\sigma \le \widecheck{(\alpha - \chi)}$ (since $\widecheck{(\alpha - \chi)}$ is the largest antichain have the non-dominating property), and hence, $\sigma \in [\chi,\widecheck{(\alpha - \chi)}]$, which establishes the set equality and completes the proof of the theorem. 
\myendproof

%Given a chain of antichains 
%$
%\{\alpha_0< \alpha_1< \alpha_2 < \dots < \alpha_k\}
%$
%such that 
%$
%\forall i \not= j \in \{0,\dots k\} : \alpha_i \cap \alpha_j = \emptyset
%$.
%The interval $[\alpha_0, \alpha_k]$ can be partitioned as
%\begin{eqnarray}
%[\alpha_0,\alpha_k] & = & \bigcup_{\chi_1 \subseteq \alpha_1} [\alpha_0 \vee \chi_1,\alpha_1^- \vee \chi_1] \ovee \bigcup_{\chi_2 \subseteq \alpha_{2,\chi_1}} [\chi_1 \vee \chi_2, \alpha_{2,\chi_1}^- \vee \chi_2] \nonumber \\
%                                    & \ldots & \ovee \bigcup_{\chi_k \subseteq \alpha_{k,\chi_{k-1}}} [\chi_{k-1} \vee \chi_k, \alpha_{k,\chi_{k-1}}^- \vee \chi_k]    
%\end{eqnarray}
%where $\forall i \in \{1,\dots k\} : \alpha_{i,\chi_{i-1}} = \{X \in \alpha_i | \{X\} \wedge \alpha_{i-1} \le \chi_{i-1}\}$.
\begin{corollary}
\label{cor:largestnondominatingdecomposition}
For any $\gamma \in [\alpha,\beta]$, we have
\[
 [\alpha,\beta] = \bigcup_{\chi \subseteq \gamma}  
 [\alpha \vee \chi,\beta \wedge \widecheck{(\gamma - \chi)}]
\]
where the intervals are pairwise disjoint.
\end{corollary}
\proof
By Theorem \ref{the:largestnondominatingdecomposition} and Proposition \ref{prop:intervalintersect}, we have
\begin{eqnarray*}
[\alpha,\beta] & = & \A_N \cap [\alpha,\beta]\\
           & = & (\bigcup_{\chi \subseteq \gamma} [\chi,\widecheck{(\gamma - \chi)}]) \cap [\alpha,\beta]\\
           & = & \bigcup_{\chi \subseteq \gamma} ([\chi,\widecheck{(\gamma - \chi)}] \cap [\alpha,\beta])\\
           & = & \bigcup_{\chi \subseteq \gamma} [\alpha\vee\chi,\beta \wedge\widecheck{(\gamma - \chi)}]
\end{eqnarray*}
\myendproof

\label{sub:compositiondirectproduct}
In \cite{PDCSDW2013, PDCSDW2011}, the direct product operator (defined in Section \ref{sec:preliminaries}) was used to decompose $\A_n$ in terms of $\A_{s}$ and $\A_{t}$ with $s + t = n$.
Here we give an alternate proof based on interval posets.
\begin{lemma}
\label{lem:projectionlemma}
Let $\chi \in \A_N$ and let $N_1 \cup N_2 = N = \{1,2,\ldots,n\}$ be a partition of $N$.
Then
\begin{equation*}
\chi \in [\alpha_1 \vee \alpha_2,\alpha_1 \varotimes \alpha_2]
\end{equation*}
where $\alpha_i = \chi \wedge \{N_i\}, i = 1,2$.
\end{lemma}
\proof
$\alpha_1 \le \chi$ and $\alpha_2 \le \chi \Rightarrow \alpha_1 \vee \alpha_2 \le \chi$.
To show $\chi \le \alpha_1 \varotimes \alpha_2$, let $A \in \chi$. We have $A = A_1 \cup A_2$, where $A_i = A \cap N_i, i=1,2$.
\begin{eqnarray*}
\{A_1\} & = & \{A \cap N_1\} \le \{A\} \le \chi       \text{ and } \{A \cap N_1\} \le \{N_1\}\\
 \Rightarrow & & \{A_1\} \le \chi \wedge \{N_1\} = \alpha_1\\
 \Rightarrow & & \exists B_1 \in \alpha_1 \text{ s.t. } A_1 \subseteq B_1
\end{eqnarray*}
By the same argument, $\exists B_2 \in \alpha_2 \text{ s.t. } A_2 \subseteq B_2$.
We have $B = B_1 \cup B_2 \in \alpha_1 \varotimes \alpha_2$, and hence, 
$A = A_1 \cup A_2 \subseteq B_1 \cup B_2 = B \in \alpha_1 \varotimes \alpha_2$,
which shows that $\chi \le \alpha_1 \varotimes \alpha_2$, and completes the proof of the lemma.
\myendproof
%Let $N = \{1,2,\ldots,n\}$ be a set of natural numbers and $\{N_1, N_2\}$ a partition of $N$.
%Given two non-empty antichains $\alpha_1 \in \A_{N_1}$ and  $\alpha_2 \in \A_{N_2}$, let $\mathcal{S}_{\alpha_1,\alpha_2} = \{A_1 \cup A_2 |  \emptyset \subsetneq A_1 \subseteq N_1, \emptyset \subsetneq  A_2 \subseteq N_2, \{A_1\} \le \alpha_1,\{A_2\} \le \alpha_2\}$. $\mathcal{S}_{\alpha_1,\alpha_2} $ is an interval poset underlying the interval $[\alpha_1 \vee \alpha_2, \alpha_1 \otimes \alpha_2]$
%\proof
%Prove that $\mathcal{S}_{\alpha_1,\alpha_2} $ satisfies the conditions of lemma \ref{the:posetcharacteristic} and note that 
%\begin{itemize}
%\item $\bigvee_{X \in \mathcal{S}_{\alpha_1,\alpha_2}} \{X\} = \alpha_1 \otimes \alpha_2$
%\item $\bigvee_{X \in \mathcal{S}_{\alpha_1,\alpha_2}} \{X^{-} \backslash \mathcal{S}_{\alpha_1,\alpha_2}\} = \alpha_1 \vee \alpha_2$
%\end{itemize}
%
%\myendproof

Note that $\alpha \otimes \beta$ is the largest antichain $\chi$ for which $\chi \wedge \{\cup{\alpha}\} = \alpha$ and $\chi \wedge \{\cup{\beta}\} = \beta$,
while $\alpha \vee \beta$ is the smallest antichain with those properties.

\begin{theorem}
\label{the:intervalintervaldecomposition}
Let $\{N_1,N_2\}$ be a partition of $N = \{1,2,\ldots,n\}$. Then
\begin{equation*}
\A_N = \{[\bot,\bot]\} \cup \bigcup_{\substack{\alpha_1 \in \A_{N_1} - \{\bot\} \\ \alpha_2 \in \A_{N_2} - \{\bot\}}} [\alpha_1 \vee \alpha_2, \alpha_1 \otimes \alpha_2] 
\end{equation*}
where the intervals are pairwise disjoint.
\end{theorem}
\proof
By Lemma \ref{lem:projectionlemma}, $\A_N$ is the union of the intervals. To show that the intervals are pairwise disjoint, let
$\chi \in [\alpha_1 \vee \alpha_2, \alpha_1 \otimes \alpha_2] \cap [\beta_1 \vee \beta_2, \beta_1 \otimes \beta_2] $, where $\chi \not= \bot$ and $\alpha_i,\beta_i \in \A_{N_i} - \{\bot\}, i = 1,2$.
The case $\chi = \bot$ is covered by noticing that $[\bot,\bot]$ is disjoint from the other intervals in the collection.
Let $A_1 \in \alpha_1$. Since $\alpha_1 \le \alpha_1 \vee \alpha_2 \le \chi$, $A \subseteq X$ for some $X \in \chi$.
But $\chi \le \beta_1 \varotimes \beta_2$ s.t. $X \subseteq B$ for some $B \in \beta_1 \varotimes \beta_2$.
Since $A_1 \subseteq N_1$, we have $A_1 \subseteq X \cap N_1 \subseteq B \cap N_1 = B_1 \in \beta_1$, which shows that $\alpha_1 \le \beta_1$.
Reversing the roles of $\alpha$ and $\beta$ in the argument above, we have $\beta_1 \le \alpha_1$, and hence, $\alpha_1 = \beta_1$. Similarly, $\alpha_2 = \beta_2$, which shows that the two intervals are the same interval, and completes the proof.
\myendproof
\section{Induced Decomposition}
\label{sec:decomposebyproperty}
In this section, the size of the sets in an antichain is used to derive decompositions and expressions for the size of (intervals in) $\A_N$.
An efficient powers-of-two formula for the size of an interval is established in Theorem \ref{the:evenoddleveldecomposition}.

Our guiding principle is the size of the sets in antichains. The counting strategy in this section is based on a certain type of decomposition in which each antichain in $[\alpha,\beta]$ is expressed as the join of $\alpha$ with the joins of {\em uniform} antichains in $\PAB$.
\begin{definition}
A {\bf uniform} antichain is an antichain whose elements (subsets of $N$) all have the same size. If that common size equals $l$, it is called {\it $l$-uniform}.
\end{definition}
Let $\PAB^l$ denote the set of all subsets in $\PAB$ of size $l$, i.e.,
\[
\PAB^l = \{X \in \PAB | |X| = l\}
\]
We will refer to $l$ as the {\it level} of the uniform antichain.
Observe that $\PP^l_{[\alpha,\beta]}$ and each of its subsets are the $l$-uniform antichains in $\PAB$.

The following two properties are immediate consequences of the definitions.
\begin{proposition}
\label{prop:luniform}
If $\gamma$ and $\delta$ are l-uniform antichains, for some $l$, then
\begin{eqnarray*}
(i)\ \gamma \le \delta & \Leftrightarrow & \gamma \subseteq \delta \\
(ii)\  \gamma \vee \delta & = & \gamma \cup \delta
\end{eqnarray*}
\end{proposition}
\begin{proposition}
\label{prop:leveldecompositionnonunique}
For two antichains $\alpha \le \beta$, any antichain $\chi \in [\alpha,\beta]$ can be decomposed as
\[
\chi = \alpha \vee \chi_m \vee \ldots \chi_M
\]
where $\chi_i \in \PP^i_{[\alpha,\beta]}$ and $m$ and $M$ are the sizes of the smallest and largest sets in $\PAB, respectively$.
\end{proposition}
% {\color{blue} Make a proposition}
\proof
By Proposition \ref{prop:posetspan}, each antichain $\chi \in [\alpha,\beta]$ can be written as $\chi = \alpha \vee \chi'$, where $\chi'$ is an antichain in $\PAB$.
We can then express $\chi'$ as the join of uniform antichains in $\PAB$. In particular, if $m$ is the size of the smallest subset in $\PAB$ and $M$
is the size of the largest, then $\chi'=\chi_m \vee \chi_{m+1} \vee \ldots \vee \chi_M$, where $\chi_i = \chi' \cap \PP^i_{[\alpha,\beta]}, i = m \ldots M$.
Thus, each antichain in $[\alpha,\beta]$ can be decomposed as $\alpha \vee \chi_m \vee \ldots \vee \chi_M$, where $\chi_i \subseteq \PP^i_{[\alpha,\beta]}$ (possibly empty), $i=m \ldots M$.
\myendproof
This decomposition is not unique, as the following example illustrates. 
\begin{example} if $\alpha = \{\{1\}\}, \beta = \{\{1,2,3\}\}$ and $\chi = \{\{1\},\{2,3\}\} \in [\alpha,\beta]$, we have
\label{ex:multipledecompositions}
\begin{eqnarray*}
\chi & = & \alpha \vee \{\{2,3\}\} \\
& = & \alpha \vee \{\{2\}\} \vee \{\{2,3\}\} \\
& = & \alpha \vee \{\{2\},\{3\}\} \vee \{\{2,3\}\}
\end{eqnarray*}
Each of these decompositions is in the form $\alpha \vee \chi_1 \vee \chi_2$, with $\chi_2=\{\{2,3\}\}$ and $\chi_1 = \{\}, \chi_1 = \{\{2\}\},\text{ or } \chi_1 = \{\{2\},\{3\}\}$ respectively.
Notice that in the last decomposition, $\chi_1$ contains both of the subsets of the one set in $\chi_2$.
If we require our decomposition to contain all such {\em redundancies}, our decomposition will be unique,
as will be shown in Theorem \ref{the:leveldecomposition}.
\end{example}
\begin{definition}
% {\color{blue} Use unions}
For an $l$-uniform antichain $\chi \subseteq \PAB$, we let \[\chi^- =\bigcup_{X \in \chi} (\pred{X} \cap \PAB)\ \text{ and }\chi^+ =\{X \in \PP^{l+1}_{[\alpha,\beta]} | \pred{X} \cap \PAB \subseteq \chi\}\]
By convention, for $[\alpha,\beta] = [\bot,\top]$, $\bot^- =  \{\emptyset\}^- = \bot$, and $\bot^+ = \{\emptyset\}$.
\end{definition}
Note that $\chi^- \le \chi$ and $(\chi^-)^+ \le \chi$. 
\begin{proposition}
\label{prop:updown}
The following properties hold for any two uniform antichains $\delta \subseteq \PP^l_{[\alpha,\beta]}$ and $\gamma \subseteq \PP^{l+1}_{[\alpha,\beta]}$.
\begin{itemize}
\item[(i)] $\gamma \subseteq \delta^+ \Leftrightarrow \gamma^- \subseteq \delta$
\item[(ii)] $(\delta^-)^+ \subseteq \delta$
\item[(iii)] $(\delta^+)^- \subseteq \delta$
\item[(iv)] $\gamma \subseteq \delta^+ \Rightarrow \gamma^+ \subseteq \delta^{++}$
\item[(v)] $\gamma^- \subseteq \delta \Rightarrow \gamma^{--} \subseteq \delta^-$
\end{itemize}
\end{proposition}
\proof
\begin{itemize}
\item[(i)] Let $\gamma \subseteq \delta^+$ and $A \in \gamma^-$. We must show $A \in \delta$.
Since  $A \in \gamma^- = \bigcup_{X \in \gamma} (\pred{X} \cap \PAB)$, we have $A \in \pred{X} \cap \PAB$ for some $X \in \gamma$.
But $\gamma \subseteq \delta^+ \Rightarrow X \in \delta^+ = \{X \in \PP^{l+1}_{[\alpha,\beta]} | \pred{X} \cap \PAB \subseteq \delta\} \Rightarrow A \in \delta$ establishing $\gamma \subseteq \delta^+ \Rightarrow \gamma^- \subseteq \delta$.

Let  $\gamma^- \subseteq \delta$ and $A \in \gamma$. We must show $A \in \delta^+$.
Since $A \in \gamma \Rightarrow \pred{A}\cap \PAB \subseteq \gamma^- \subseteq \delta$.
Since $\delta^+ =  \{X \in \PP^{l+1} | \pred{X} \cap \PAB \subseteq \delta\}$ it follows $A \in \delta^+$ establishing $\gamma \subseteq \delta^+ \Leftarrow \gamma^- \subseteq \delta$ and completing the proof of (i).
\item[(ii)] Suppose that $A \in (\gamma^-)^+$. Then $A \in \{X \in \PP^l_{[\alpha,\beta]} | \pred{X} \cap \PAB \subseteq \gamma^-\}$.
Then $A \in \PP^l_{[\alpha,\beta]}$ and $\pred{A} \cap \PAB \subseteq \gamma^-$.
But $\gamma^- = \bigcup_{X \in \gamma} (\pred{X} \cap \PAB)$.
If $A \not\in \gamma$, then $\pred{A} \cap \PAB \not\subseteq \gamma^-$, a contradiction.
Thus $A \in \gamma$, completing the proof of (ii). 
\item[(iii)] Let $\rho = \delta^+$. Then $\rho \subseteq \delta^+$ and by (i), $\rho^- \subseteq \delta$. Thus, $(\delta^+)^- \subseteq \delta$.
\item[(iv)] Suppose $\gamma \subseteq \delta^+$ and let $A \in \gamma^+$. By definition of the $()^+$ operator, $A \in \{X \in \PP^{l+2}_{[\alpha,\beta]} | \pred{X} \cap \PAB \subseteq \gamma \}$, and hence $A \in \{X \in \PP^{l+2}_{[\alpha,\beta]} | \pred{X} \cap \PAB \subseteq \delta^+ \} = \delta^{++}$. Thus, $\gamma^+ \subseteq \delta^{++}$.
\item[(v)] This follows from (i) and (iv).
\end{itemize}
\myendproof

% {\color{blue} Corollary ++++ and -----}
\begin{corollary}
For any two uniform antichains $\delta \subseteq \PP^l_{[\alpha,\beta]}$ and $\gamma \subseteq \PP^{l+j}_{[\alpha,\beta]}$:
\begin{eqnarray*}
(i)\ \gamma \subseteq \delta^+ & \Rightarrow & \gamma^{+j} \subseteq \delta^{+(j+1)} \\
(ii)\ \gamma^- \subseteq \delta & \Rightarrow & \gamma^{+(j-1)} \subseteq \delta^{+j}
\end{eqnarray*}
\end{corollary}
\proof
Assertion $(i)$ follows by repeated applications of Proposition \ref{prop:updown}$(iv)$, and Assertion $(ii)$ follows from Proposition \ref{prop:updown}$(i)$ and Assertion $(i)$.
\myendproof

As we observed in the previous example, the non-uniqueness of the decomposition resulted from being able to vary the redundancies appearing in one or more of the $\chi_i$'s.
The following theorem establishes uniqueness by imposing the condition $\chi_{i+1}^- \subseteq \chi_i$, which forces all possible redundancies with respect to $\chi_{i+1}$ to appear in $\chi_i$.

For the remainder of this paper, m and M will denote the size of the smallest and the largest subsets in $\PAB$.
\begin{theorem}
\label{the:leveldecomposition}
For each $\chi \in [\alpha,\beta]$, there is exactly one decomposition of the form $\chi = \alpha \vee \chi_m \vee \ldots \vee \chi_M$, where $\chi_{i+1}^- \subseteq \chi_i \subseteq \PP^i_{[\alpha,\beta]}$, for $i = m, \ldots, M-1$.
% \begin{itemize}
% \item for $\chi_i \not= \bot$, we have $\chi_i \subseteq \PP^i_{[\alpha,\beta]}$ and $\chi_{i+1} \subseteq \chi_i^+$
% \item for $\chi_i = \bot$, $\chi_{i+1} \subseteq \{X \in \PP^{i+1}_{[\alpha,\beta]} | \pred{X} \cap \PAB = \emptyset\}$
% \end{itemize}
\end{theorem}
\proof
Let $\chi \in [\alpha,\beta]$ and consider $\alpha \vee \chi_m \vee \ldots \vee \chi_M$,  where  \hfill\break {$\chi_i = [(\chi - \alpha) \cap \PP^i_{[\alpha,\beta]}] \cup \chi_{i+1}^-$, for $i = m, \ldots, M-1$}.
The definition of $\chi_i$ implies $\chi = \alpha \vee \chi_m \vee \ldots \vee \chi_M$ and $\chi_{i+1}^- \subseteq \chi_i \subseteq \PP^i_{[\alpha,\beta]}$.
It remains to show that this decomposition is unique.

Suppose $\alpha \vee \chi_m \vee \ldots \vee \chi_M = \alpha \vee \gamma_m \vee \ldots \vee \gamma_M$ are two decompositions of the antichain $\chi$,
where $\chi_i \subseteq  \PP^i_{[\alpha,\beta]}$ and $\gamma_i \subseteq  \PP^i_{[\alpha,\beta]}$, $i=m,\ldots,M$, and $\chi_{i+1}^- \subseteq \chi_i$ and $\gamma^-_{i+1} \subseteq \gamma_i$, $i = m,\dots,M-1$.
We must show $\chi_i = \gamma_i$, $i = m,\ldots, M$.
Since M is the size of the largest subset appearing in any antichain in $[\alpha,\beta]$, every subset appearing in $\chi_M$ must appear in $\gamma_M$ and vice versa,
otherwise, the two decompositions would not be equal. Thus, $\chi_M = \gamma_M$.
To show that $\chi_{M-1} = \gamma_{M-1}$,
let $A \in \chi_{M-1}$. Since $\alpha \vee \chi_m \vee \ldots \vee \chi_M = \alpha \vee \gamma_m \vee \ldots \vee \gamma_M$, $A \subseteq B$ for some $B \in \alpha \vee \gamma_m \vee \ldots \vee \gamma_M$.
If $A = B$, then $A \in \gamma_{M-1}$.
If $A \subsetneq B$, then $|B| = M$ and $B \in \gamma_M$. Thus $A \in pred(B) \subseteq \gamma^-_M \subseteq \gamma_{M-1}$, which shows that $\chi_{M-1} \subseteq \gamma_{M-1}$.
A symmetric argument shows the reverse subset relation, and hence, $\chi_{M-1} = \gamma_{M-1}$.

Continuing in this way, one can show $\chi_i = \gamma_i$ for all $i$, thereby establishing uniqueness.
%Let $\chi_M$ be the M\_uniform antichain consisting of all sets of $\chi - \alpha$ of size $M$.
%For $i = m+1 \ldots M$, given $\chi_i$, let $\chi_{i-1}$ be the (i-1)\_uniform antichain consisting of the sets in $\chi_i^-$ and all sets of $\chi-\alpha$ of size $i-1$.
%Then $\chi = \alpha \vee \chi_m \vee \ldots \chi_M$.
%If $\ch_ii^- \not= \bot$, then since $(\chi^-_i)^+ \subseteq \chi_i$, the conditions are satisfied.
%if $\chi_i^-  = \bot$, it is possible that $\chi_{i-1} = \bot$. The second condition then is satisfied.  
%So there is an expansion satisfying the conditions of the proposition.
%
%Given the expansion, we construct the $\chi^{ }_i$ recursively going down from $i=M$ to $i=m$ as follows.
%Consider a set of size $M$ in $\chi - \alpha$. Such a set must occur in $\chi_M$.
%Conversely, given a set in $\chi^{ }_M$, this set is not in $\alpha$ and must occur in $\chi$.
%We hence have determined $\chi_M$ as the collection of sets of size $M$ in $\chi - \alpha$.
%Suppose we have proven uniqueness of all $\chi_i$ for $i > l \ge m$.
%In case $\chi_{l+1}^- \not= \bot$, due to the condition $\chi_{l+1} \subseteq \chi_l^+$, we have $\chi_{l+1}^- \le \chi_l$ and since $\chi_l$ and $\chi_{l+1}$ are l(l+1)\_uniform antichains respectively, $\chi_{l+1}^- \subseteq \chi_l$.
%In this case as well as in case  $\chi_{l+1}^- \not= \bot$, any set in $\chi_l$ that does not occur in $\chi^-_{l+1}$ wil occur in $\chi-\alpha$ establishing uniqueness .
\myendproof

The following example illustrates how we can use Theorem \ref{the:leveldecomposition} to calculate the size of an interval.
Observe that the condition
$
\chi_{i+1} \subseteq \chi_i^+$ implies
$
\chi_{i+1}^+ \subseteq \chi_i^{++}
$
by Proposition \ref{prop:updown}$(iv)$.

\begin{example}
Let $\alpha = \{\{1\}\}$ and $\beta = \{\{1,2,3\}\}$. Then
\[
\PAB = \{\{2\},\{3\},\{1,2\},\{1,3\},\{2,3\},\{1,2,3\}\}\text{, with } m = 1\text{ and }M = 3.
\]
By Theorem \ref{the:leveldecomposition}, $|[\alpha,\beta]|$ equals the number of ways we can form the join $\alpha \vee \chi_1 \vee \chi_2 \vee \chi_3$ such that 
$\chi_i \subseteq \PAB$ and $\chi_{i+1}^- \subseteq \chi_i$ for $i = 1,2$.
Note that $\chi_{i+1}^- \subseteq \chi_i \Rightarrow \chi_{i+1}^{--} \subseteq \chi_{i}^-$ so that we have $\chi_3^{--} \subseteq \chi_1$.
There are two choices for $\chi_3$: $\chi_3 = \{\}$ or $\chi_3 = \{\{1,2,3\}\}$.
If $\chi_3 = \{\}$, then we have $\chi_3^- = \chi_3^{--} = \{\}$, which does not impose any condition on $\chi_1$ and $\chi_2$. There are four possible choices for $\chi_1$, and for each, we have the following possibilities for $\chi_2$:
\begin{itemize}
\item {$\chi_1 = \{\}$, $\chi_2^- \subseteq \chi_1 \Rightarrow \chi_2 = \{\}$, 1 possibility}
\item {$\chi_1 = \{\{2\}\}$, $\chi_2^- \subseteq \chi_1 \Rightarrow \chi_2 \subseteq \{\{1,2\}\}$, 2 possibilities}
\item {$\chi_1 = \{\{3\}\}$, $\chi_2^- \subseteq \chi_1 \Rightarrow \chi_2 \subseteq \{\{1,3\}\}$, 2 possibilities}
\item {$\chi_1 = \{\{2\},\{3\}\}$, $\chi_2^- \subseteq \chi_1 \Rightarrow \chi_2 \subseteq \{\{1,3\},\{1,2\},\{2,3\}\}$, 8 possibiillties}
\end{itemize}

If $\chi_3 = \{\{1,2,3\}\}$, then $\chi_3^- = \{\{1,2\},\{1,3\},\{2,3\}\}, \chi_3^{--} = \{\{2\},\{3\}\}$, which forces $\chi_1 = \{\{2\},\{3\}\}$ and 
$\chi_2 = \{\{1,2\},\{1,3\},\{2,3\}\}$.
Thus, we have a total of $1 + 2 + 2 + 8 + 1$ possibilities, and hence, $|[\alpha,\beta]| = 14$.
\end{example}
% {\color{blue} put lemma 5.8 and theorem 5.9 here}
The next result is a key element for the main counting result of this section.
\begin{lemma}
\label{lem:powerlevel}
In the construction of a decomposition $\alpha \vee \chi_m \vee \ldots \vee \chi_M$, if the choices for $\chi_i$ and $\chi_{i+2}$ have been made, then the number of possible choices for $\chi_{i+1}$ equals $2^{|\chi_i^+| - |\chi_{i+2}^-|}$.
\end{lemma}
\proof
If $\chi_i$ is fixed, the condition $\chi_{i+1}^- \subseteq \chi_i$ implies (by Proposition \ref{prop:updown}(i)) $\chi_{i+1} \subseteq \chi_i^+$.
With the condition $ \chi_{i+2}^- \subseteq \chi_{i+1}$ we arrive at $\chi_{i+2}^- \subseteq \chi_{i+1} \subseteq \chi_i^+$ leaving exactly $2^{|\chi_i^+| - |\chi_{i+2}^-|}$ possibilities for $\chi_{i+1}$.
\myendproof

%\begin{definition}
%For $\chi_i \subseteq \PP_{[\alpha,\beta]}^i$, $\chi_{i+2} \subseteq \PP_{[\alpha,\beta]}^{i+2}$ and $\chi_{i+2}^{--} \subseteq \chi_i$, the function $\bet$ is defined as
%\[
%\bet(\chi_i,\chi_{i+2}) = 2^{|\chi_i^+| - |\chi_{i+2}^-|}
%\]
%\end{definition}
%The value of $\bet(\chi_i,\chi_{i+2})$ is the number of antichains $\chi_{i+1} \subseteq \PP_{[\alpha,\beta]}^{i+1}$ for which $\chi_{i+2}^- \subseteq \chi_{i+1} \subseteq \chi_i^+$.
\begin{theorem}
\label{the:evenoddleveldecomposition}
Let $r = \lfloor \frac{M-m}{2} \rfloor$ and $\Delta = (M-m)\ mod\ 2$. We have the following two ways of calculating $|[\alpha,\beta]|$$:$
\begin{flalign*}
%(i)\ |[\alpha,\beta]| & =  \\
 % \sum_{\chi_0 \subseteq \PAB^M} & \sum_{\chi_{0}^{--} \subseteq \chi_{1}} \ldots \sum_{\chi_{r-1}^{--} \subseteq \chi_{r}}  
% (\prod_{i=0}^{r-1}\bet(\chi_{i+1},\chi_{i}))* 2^{\Delta*(|\PAB^m| - |\chi_{m+1}^-|)}\\
(i)\ |[\alpha,\beta]| & =  \\
  \sum_{\chi_M \subseteq \PAB^M} & \sum_{\chi_{M}^{--} \subseteq \chi_{M-2}} \ldots \sum_{\chi_{M-2r+2}^{--} \subseteq \chi_{M-2r}}  
% (\prod_{i=0}^{r-1}\bet(\chi_{M-2i-2},\chi_{M-2i}))* 2^{\Delta*(|\PAB^m| - |\chi_{m+1}^-|)}\\
 2^{\sum_{i=0}^{r-1}(|\chi_{M-2i-2}^+| - |\chi_{M-2i}^-|) + \Delta*(|\PAB^m| - |\chi_{m+1}^-|)} \\
%(i)\ |[\alpha,\beta]| & =  \\
%  \sum_{\chi_M \subseteq \PAB^M} & \sum_{\chi_{M}^{--} \subseteq \chi_{M-2}} 2^{|\chi_{M-2}^+| - |\chi_{M}^-|} 
 % \sum_{\chi_{M-2}^{--} \subseteq \chi_{M-4}} 2^{|\chi_{M-4}^+| - |\chi_{M-2}^-|}\\
%  & \ldots \sum_{\chi_{M-2r+2}^{--} \subseteq \chi_{M-2r}}  2^{|\chi_{M-2r}^+| - |\chi_{M-2r+2}^-|) + \Delta*(|\PAB^m| - |\chi_{m+1}^-|)}\\
% (\prod_{i=0}^{r-1}\bet(\chi_{M-2i-2},\chi_{M-2i}))* 2^{\Delta*(|\PAB^m| - |\chi_{m+1}^-|)}\\
%2^{\sum_{i=0}^{r-1}(|\chi_{M-2i-2}^+| - |\chi_{M-2i}^-|) + \Delta*(|\PAB^m| - |\chi_{m+1}^-|)} \\
(ii)\ |[\alpha,\beta]| & =  \\
\sum_{\chi_{M-1} \subseteq \PAB^{M-1}} & \sum_{\chi_{M-1}^{--} \subseteq \chi_{M-3}} \ldots \sum_{\chi_{M-2r+3}^{--} \subseteq \chi_{M-2r+1}} \\
% & \ \ \ \ \ \ \ 2^{|\chi_{M-1}^+|} * (\prod_{i=1}^{r-1} \Theta(\chi_{M-1-2i}^+,\chi_{M + 1 - 2i}^-)) * 2^{(1-\Delta)*(|\PAB^m| - |\chi_{m+1}^-|)} \\
 & \ \ \ \ \ \ \ 2^{|\chi_{M-1}^+| + | \sum_{i=1}^{r-1}(|\chi_{M-1-2i}^+| - |\chi_{M + 1 - 2i}^-|) + (1-\Delta)*(|\PAB^m| - |\chi_{m+1}^-|)} \\
\end{flalign*}
\end{theorem}
\proof
$(i)$ For each choice of $\chi_{M-2i}$, $i=0,1,\ldots,r-1$, the number of choices for $\chi_{M-2i-1}$ equals $2^{|\chi_{M-2i-2}^+| - |\chi_{M-2i}^-|}$ by Lemma \ref{lem:powerlevel}.
It follows that the total number of choices for any such $\chi_{M-2i-1}$ equals 
\[
2^{\sum_{i=0}^{r-1} |\chi_{M-2i-2}^+| - |\chi_{M-2i}^-|}
\]
The quantity $\Delta$ accounts for the extra term that occurs when $M-m$ is odd.
Equation $(ii)$ follows by the same argument.
\myendproof

Theorem \ref{the:evenoddleveldecomposition} is illustrated in Algorithms \ref{alg:evenoddleveldecomposition} and \ref{alg:evenoddleveldecompositionbasic}.
It enables us to compute the size of any interval efficiently by summing over half of the uniform-antichain levels only. 
Moreover, the two formulas provide flexibility in approaching various problems. 
For example, when computing the size of $\A_6$, the size of $\PP_{\A_6}^3$ is $20$ with $2^{20}$ subsets. 
Formula \ref{the:evenoddleveldecomposition}(i) sums over subsets of $\PP_{\A_6}^6, \PP_{\A_6}^4,\PP_{\A_6}^2$ and $\PP_{\A_6}^0$, avoiding the largest level $\PP_{\A_6}^3$. 
In this case, the range-sizes for the sums in formula \ref{the:evenoddleveldecomposition}(i) are always less than $2^{\binom{6}{2}} = 2^{15}$. In general, when computing the size of an interval, avoiding the largest level leads to a more efficient computation and motivates the choice between \ref{the:evenoddleveldecomposition}(i) or \ref{the:evenoddleveldecomposition}(ii). 

\begin{algorithm}
\caption {Compute the size of an interval $[\alpha,\beta]$ following Theorem \ref{the:evenoddleveldecomposition}}
\label{alg:evenoddleveldecomposition}
\begin{algorithmic}
\IF {(Algorithm  \ref{the:evenoddleveldecomposition}$(i)$ is the most efficient option)}
\STATE \COMMENT{Compute the size using Theorem \ref{the:evenoddleveldecomposition}$(i).$}
\STATE{$Sum \gets 0$}
\FOR{$\chi_M \subseteq \PAB^M$}
\STATE{$Sum \gets Sum + NumberAntiChains(\chi_M)$}
\STATE \COMMENT{$NumberAntiChains(\chi_M)$ is the number of antichains $\gamma$ such that the sizes of the sets in $\gamma$ are $\le M$ and $\gamma \cap \PAB^{M}== \chi_{M}$.}
\ENDFOR
\RETURN{Sum}
\ELSE
\STATE \COMMENT{Compute the size using Theorem \ref{the:evenoddleveldecomposition}$(ii)$.}
\STATE{$Sum \gets 0$}
\FOR{$\chi_{M-1} \subseteq \PAB^{M-1}$}
\STATE{$Sum \gets Sum + 2^{|\chi_{M-1}^+|}*NumberAntiChains(\chi_{M-1})$}
\ENDFOR
\RETURN{Sum}
\ENDIF
\end{algorithmic}
\end{algorithm}

\begin{algorithm}
\caption {NumberAntiChains($\chi_{l}$)}
\label{alg:evenoddleveldecompositionbasic}
\begin{algorithmic}
\STATE {\bf Input} {$m \le l \le M;\ \chi_{l} \subseteq \PAB^{l}$}
%\STATE \COMMENT{Return the number of antichains $\gamma \in [\alpha,\beta]$ such that the sizes of the sets in $\gamma$ are $\le l$ and $\gamma \cap \PAB^{l}== \chi_{l}$ }
\IF{$l == m$}
\RETURN{$1$}
\ENDIF
\IF{$l == m+1$}
\RETURN{$2^{|\PAB^m| - |\chi_{l}^-|}$}
\ENDIF
\STATE{$Sum \gets 0$}
\FOR {$\chi_{l-2} \supseteq \chi_{l}^{--}$}
	\STATE{$Sum \gets Sum + 2^{|\chi_{l-2}^+| - |\chi_{l}^-|}*NumberAntiChains(\chi_{l-2})$}
\ENDFOR
\RETURN $Sum$
\end{algorithmic}
\end{algorithm}

\begin{samepage}
\begin{lemma}
\label{lem:subset}
Let  $\chi = \alpha \vee \chi_m \vee \ldots \vee \chi_M$ be the decomposition of $\chi \in [\alpha,\beta]$, specified by Theorem \ref{the:leveldecomposition}.
If $\{A\} \le \chi_i$ for some $i$, and $|A| = l$, then $\{A\} \le \alpha$ or $A \in \chi_l$.
\end{lemma}
\proof
Since $\{A\} \le \chi_i$, we have $\alpha \vee \{A\} \in [\alpha,\beta]$, and by definition of $\PAB$, either $\{A\} \le \alpha$ (in which case we are done)
or $A \in \PAB$.
To complete the proof for the latter case, we must show that $A \in \chi_l$.
We have $A \subseteq B$ for some $B \in \chi_i$ (since $\{A\} \le \chi_i$).
If $A = B$, then $i = l$, and we are done.
Otherwise, $A = B - \{b_1,b_2,\dots,b_{i-l}\}$ where $\{b_1,b_2,\dots,b_{i-l}\} \subseteq B$.
Thus, $A \subseteq B - \{b_1,b_2,\dots,b_{k}\} \subseteq B$ for $k = 1,\ldots,i-l$, where $A \in \PAB$ and also $B \in \PAB$ (since $B \in \chi_i$).
By Theorem \ref{the:posetcharacteristic}, the subsets $B - \{b_1\}, B - \{b_1,b_2\},\ldots,B-\{b_1,b_2,\ldots,b_{i-l}\}$ are all elements in $\PAB$,
and by definition of the $()^-$ operator, $B-\{b_1\} \in \chi_i^-$. But $\chi_i^- \subseteq \chi_{i-1}$, and hence, $B-\{b_1\} \in \chi_{i-1}$. 
Similarly, $B-\{b_1,b_2,\ldots,b_k\} \in \chi_{i-k}$ for $k=1,\ldots,i-l$. In particular, $A = B-\{b_1,b_2,\dots,b_{i-l}\} \in \chi_l$,
which concludes the proof.
\myendproof
\end{samepage}
\pagebreak
\begin{proposition}
\label{prop:wedgeintersection}
Let $\chi = \alpha \vee \chi_m \vee \ldots \chi_M$ and $\gamma = \alpha \vee \gamma_m \vee \ldots \gamma_M$ be the decompositions specified by Theorem \ref{the:leveldecomposition} of two antichains $\chi,\gamma \in [\alpha,\beta]$. Then
\[
\chi \wedge \gamma = \alpha \vee (\chi_m \cap \gamma_m) \vee \dots \vee (\chi_M \cap \gamma_M)
\]
\end{proposition}
\proof
Using the distributive property and the inequalities $\alpha \le \chi$ and $\alpha \le \gamma$, we have
\[
\chi \wedge \gamma = \alpha \vee \bigvee_{i,j \in \{m,\dots,M\}} (\chi_i \wedge \gamma_j)
\]
To complete the proof, we must show
\[
 \alpha \vee \bigvee_{i,j \in \{m,\dots,M\}} (\chi_i \wedge \gamma_j) =  \alpha \vee \bigvee_{l \in \{m,\ldots,M\}} (\chi_l \cap \gamma_l)
\]
The inequality $\ge$ is immediate from the definitions.
To establish the reverse inequality, it suffices to show that for all $i \not= j$, $\chi_i \wedge \gamma_j \le \alpha \vee (\chi_l \cap \gamma_l)$ for some $l, m \le l \le M$.
Without loss of generality, we may assume $i < j$.

Suppose $X \in \chi_i \wedge \gamma_j$.
By Proposition \ref{prop:max}, $X = A \cap B$ for some $A \in \chi_i$ and $B \in \gamma_j$, with $|A| = i < j = |B|$.
It follows that  $\{A \cap B\} \le \{A\} \le \chi_i$ and $\{A \cap B\} \le \{B\} \le \gamma_j$.
Thus, by Lemma \ref{lem:subset}, $\{A \cap B\} \le \alpha$ or $A \cap B \in \chi_l \cap \gamma_l$, which shows that $\chi_i \wedge \gamma_j \le \alpha \vee (\chi_l \cap \gamma_l)$, 
and completes the proof.
\myendproof

\begin{corollary}
\label{cor:lessthandecomposition}
Let $\chi = \alpha \vee \chi_m \vee \ldots \chi_M$ and $\gamma = \alpha \vee \gamma_m \vee \ldots \gamma_M$ be the unique decompositions specified by Theorem \ref{the:leveldecomposition} of two antichains $\chi,\gamma \in [\alpha,\beta]$. If $\chi \le \gamma$, then $\chi_i \subseteq \gamma_i$ for all $i = m,\dots,M$.
\end{corollary}
\proof
If $\chi \le \gamma$, then $\chi = \chi \wedge \gamma$, and by Proposition \ref{prop:wedgeintersection},
\[
\chi = \alpha \vee (\chi_m \cap \gamma_m) \vee \ldots \vee (\chi_M \cap \gamma_m)
\]
Theorem \ref{the:leveldecomposition} implies $\chi_i = \chi_i \cap \gamma_i, i = m,\ldots,M$, and hence, $\chi_i \subseteq \gamma_i$
% {\color{blue} See photo 5.7}
%Suppose $A \in \chi_i$ for some $i \in \{m,\dots,M\}$. We must show that $A \in \gamma_i$.
%Since $\chi_i \le \chi \le \gamma$, $A \subseteq B$ for some $B \in \gamma$. Thus, $|A| = i, |B| = j$ and $B \in \gamma_j$ for some $j \ge i$.
%By Lemma \ref{lem:subset}, $\{A\} \le \alpha$ or $A \in \gamma_i$. Since $A \in \chi_i \subseteq \PAB$, $\{A\} \not\le \alpha$, and hence, $A \in \gamma_i$.
\myendproof

%\begin{notation}
%To simplify the notation for the formulas to follow, we will use $\PAB^l$ instead of $\PP^l_{[\alpha,\beta]}$ when the interval under consideration is clear.
%\end{notation}
Theorem \ref{the:evenoddleveldecomposition} shows how to sum over all possible choices for $\chi_m,\chi_{m+2},\ldots$ given a set of choices for $\chi_{m+1},\chi_{m+3},\ldots$.
One  can consider other sets of such choices. For example, if $m < k <M$, we could consider all possible choices for $\chi_{k}$, and for each such choice compute the number of antichains $\chi \in [\alpha,\beta]$.

% {\color{blue} combine both lemma into a proposition}
We now prove a proposition that leads to a generalization of Theorem \ref{the:evenoddleveldecomposition}.
\begin{proposition}
\label{prop:antichainsinterval}
% For any $\chi \in [\alpha,\beta]$, let $\chi = \alpha \vee \chi_m \vee \ldots \vee \chi_M$ with $\chi^-_{i+1} \subseteq \chi_i \subseteq \PAB^i$ be the decomposition specified in Theorem \ref{the:leveldecomposition}. If $\rho_k \subseteq \PAB^k$ for some $k, m < k < M$, then
%\[
%\chi_k = \rho_k \Leftrightarrow \chi \in [\alpha \vee \rho_k, \alpha \vee \PAB^m \vee \ldots \vee \PAB^{k-1} \vee \rho_k \vee \rho_k^+ \vee \ldots \vee \rho_k^{+(M-k)}] 
%\]
Let $\chi$ be any antichain in the interval $I=[\alpha,\beta]$. For each $k$, $m < k < M$, $\chi_k$ is the $k$-uniform antichain in the decomposition of $\chi$ specified by Theorem \ref{the:leveldecomposition} if and only if
\[
\chi \in  [\alpha \vee \chi_k, \alpha \vee \PP_I^m \vee \ldots \vee \PP_I^{k-1} \vee \chi_k \vee \chi_k^+ \vee \ldots \vee \chi_k^{+(M-k)}]
\]
\end{proposition}
\proof

($\Rightarrow$) 
Suppose $\chi = \alpha \vee \chi_m \vee \ldots \vee \chi_M$ is the decomposition specified by Theorem \ref{the:leveldecomposition}, and consider any $k$, $m<k<M$.
The inequality $\alpha \vee \chi_k \le \chi$ is immediate.
The inequality 
\[
\alpha \vee \chi_m \vee \ldots \vee \chi_{k-1} \vee \chi_k \vee \chi_{k+1} \vee \ldots \chi_M \le
\alpha \vee \PP_I^m \vee \ldots \vee \PP_I^{k-1} \vee \chi_k \vee \chi_k^+ \vee \ldots \vee \chi_k^{+(M-k)}
\]
follows directly from the condition $\chi^-_{i+1} \subseteq \chi_i \subseteq \PP^i_I$, $i=m,\ldots,M-1$, Proposition \ref{prop:updown}$(i)$, and repeated applications of Proposition \ref{prop:updown}$(iv)$.

($\Leftarrow$) 
Suppose $\chi \in  [\alpha \vee \rho_k, \alpha \vee \PP_I^m \vee \ldots \vee \PP_I^{k-1} \vee \rho_k \vee \rho_k^+ \vee \ldots \vee \rho_k^{+(M-k)}]$, where $\rho_k \subseteq \PP_I^k$ for some $k$, $m<k<M$, and let $\chi = \alpha \vee \chi_m \vee \ldots \vee \chi_M$ be its decomposition specified by Theorem \ref{the:leveldecomposition}.
We must show $\chi_k = \rho_k$. 
First, we have
\[
\alpha \vee \chi_m \vee \ldots \vee \chi_k \vee \chi_{k+1} \vee \ldots \vee \chi_M 
= \chi \le \alpha \vee \PP_I^m \vee \ldots \vee \PP_I^{k-1} \vee \rho_k \vee \rho_k^+ \vee \ldots \vee \rho_k^{+(M-k)} 
\]
and hence, by Corollary \ref{cor:lessthandecomposition}, $\chi_k \subseteq \rho_k$.
To establish the reverse inclusion, consider the antichain $\alpha \vee \rho_k$. Using its decomposition specified by 
Theorem \ref{the:leveldecomposition}, we have
\[
\alpha \vee \rho_k^{-(k-m)} \vee \ldots \vee \rho_k^- \vee \rho_k = \alpha \vee \rho_k 
\le \chi = \alpha \vee \chi_m \vee \ldots \vee \chi_k \vee \ldots \vee \chi_M
\]
Thus, $\rho_k \subseteq \chi_k$ by Corollary \ref{cor:lessthandecomposition}.
\myendproof
\begin{lemma}
\label{lem:splitwedge}
If $\chi_k \subseteq \PAB^k$, then
\[
(\alpha \vee \PAB^m \vee \dots \PAB^{k-1} \vee \chi_k) \wedge (\alpha \vee \chi_k \vee \chi_k^+ \vee \dots \vee \chi_k^{+(M-k)}) = \alpha \vee \chi_k
\]
\end{lemma}
\proof
By Theorem \ref{the:leveldecomposition} and the definitions of the $()^-$ and $()^+$ operators, we have 
\[
\alpha \vee \chi_k \vee \chi_k^+ \vee \dots \vee \chi_k^{+(M-k)} = \alpha \vee \chi_k^{-(k-m)} \vee \dots \vee \chi_k^- \vee \chi_k \vee \chi_k^+ \vee \dots \vee \chi_k^{+(M-k)}
\]
and hence, by Proposition \ref{prop:wedgeintersection},
\begin{eqnarray*}
(\alpha \vee \PAB^m \vee \dots &\vee& \PAB^{k-1} \vee \chi_k)  \wedge (\alpha \vee \chi_k \vee \chi_k^+ \vee \dots \vee \chi_k^{+(M-k)}) \\
& = & \alpha \vee (\PAB^m \cap \chi_k^{-(k-m)}) \vee \dots \vee (\PAB^{k-1} \cap \chi_k^-) \vee (\chi_k \cap \chi_k) \\
& = & \alpha \vee \chi_k^{-(k-m)} \vee \dots \vee \chi_k^- \vee \chi_k \\
& = & \alpha \vee \chi_k 
\end{eqnarray*}
\myendproof

\begin{theorem}
\label{the:generalleveldecomposition}
Let $m < k < M$. We have
\[
|[\alpha,\beta]| = \sum_{\chi_k \subseteq \PAB^k} |[\alpha \vee \chi_k,\alpha \vee \PAB^m \vee \ldots \PAB^{k-1} \vee \chi_k ]|.|[\alpha \vee \chi_k,\alpha \vee \chi_k \vee \chi_k^+ \vee \ldots \vee \chi_k^{+(M-k)}]|  
\]
\end{theorem}
\proof
By Proposition \ref{prop:antichainsinterval}, the set of antichains in $[\alpha,\beta]$ for which $\chi_k = \rho_k \subseteq \PAB^k$ is given by
\begin{equation}
\label{eq:rhoatom}
[\alpha \vee \rho_k,\alpha \vee \PAB^m \vee \ldots \PAB^{k-1} \vee \rho_k \vee \rho_k^+ \vee \ldots \vee \rho_k^{+(M-k)}].
\end{equation}
Moreover, these intervals form a partition of $[\alpha,\beta]$, and we have
\[
|[\alpha,\beta]| = \sum_{\chi_k \subseteq \PAB^k} |[\alpha \vee \chi_k,\alpha \vee \PAB^m \vee \ldots  \PAB^{k-1} \vee \chi_k \vee \chi_k^+ \ldots \vee \chi_k^{+(M-k)}]|.
\]
Lemma \ref{lem:splitwedge} and Theorem \ref{the:directjoindecomposition} imply
% {\color{blue} Align}
\begin{multline*}
[\alpha \vee \chi_k,\alpha \vee \PAB^m \vee \ldots  \PAB^{k-1} \vee \chi_k \vee \chi_k^+ \ldots \vee \chi_k^{+(M-k)}]\\
= [\alpha \vee \chi_k,\alpha \vee \PAB^m \vee \ldots  \PAB^{k-1} \vee \chi_k]
 \varovee  [\alpha \vee \chi_k,\alpha \vee \chi_k \vee \chi_k^+ \ldots \vee \chi_k^{+(M-k)}]
\end{multline*}
and the result follows. 
\myendproof

% {\color{blue} put theorem 5.13 here}
\begin{lemma}
\label{lem:lambdaischik}
Let $I = [\alpha \vee \chi_k,\alpha \vee \chi_k \vee \chi_k^+ \ldots \vee \chi_k^{+(M-k)}]$, where $\chi_k \subseteq \PAB^k$, for some $k \in \{m+1,\ldots,M-1\}$. 
Then $\PP_I^l = \chi_k^{+(l-k)}$, for $l = k+1,\ldots,M$.
\end{lemma}
\proof
Suppose $X \in \chi_k^{+(l-k)}$. 
Then $|X| = k + l - k = l$ and $\{X\} \le \alpha \vee \chi_k \vee \chi_k^+ \vee \ldots \vee \chi_k^{+(l-k)} \vee \ldots \vee \chi_k^{+(M-k)}$.
Moreover, since $X \in \PAB$ and $|X| > k$ we have  $\alpha \vee \{X\} > \alpha$ and $\chi_k \vee \{X\} > \chi_k$, and hence, $\alpha \vee \chi_k \vee \{X\} > \alpha \vee \chi_k$.
Thus, $X \in \PP_I^l$, which shows that $\chi_k^{+(l-k)} \subseteq \PP_I^l$.

Now suppose $X \in \PP_I^l$. 
By definition, we have $\{X\} \not\le \alpha \vee \chi_k$ and $\{X\} \le \alpha \vee \chi_k  \vee \chi_k^+ \vee \ldots \vee \chi_k^{+(M-k)}$.
Moreover, since $|X| = l$, $\{X\} \not\le  \chi_k^{+(i-k)}$ for all $i < l$.
It follows that $\{X\} \le \chi_k^{+(l-k)} \vee \ldots \vee \chi^{+(M-k)}$,
and hence,  $X \subseteq Y$ for some $Y \in \chi_k^{+(i-k)}$, where $l \le i \le M$.
Among all such $Y's$, let $Y^*$ be a minimal one.
We have $X,Y^* \in \PP_I$.
If $X \subsetneq Y^*$, then $X \subseteq Y^*-\{y\}$ for some $y \in Y^*$. 
Since $X, Y^* \in \PP_I$, Theorem \ref{the:posetcharacteristic} implies that $Y^*-\{y\} \in \PP_I$, 
and by definition of the $()^-$ operator and Proposition \ref{prop:updown}$(iii)$,
%and by definition of the $()^+$ operator, 
$Y^*-\{y\} \in \chi_k^{+(i-k)-} = \chi_k^{+(i-1-k)}$,
which would contradict the minimality of $Y^*$.
Therefore, $X = Y^*$, which shows that $X \in \chi_k^{+(l-k)}$ and completes the proof.
\myendproof

\begin{notation}
Let $k_1,\ldots,k_t$ be any t integers such that $m < k_1 < k_2 < \ldots < k_t < M$, and let $\chi_{k_j} \subseteq \PAB^{k_j}$, 
$j = 1,\ldots,t$. 
\begin{itemize}
\item[$(i)$] $B_{k_1}$ denotes the interval $[\alpha \vee \chi_{k_1}, \alpha \vee \PAB^m \vee \ldots \vee \PAB^{k_1 - 1} \vee \chi_{k_1}]$.
\item[$(ii)$] $I_{k_1} = [\alpha \vee \chi_{k_1},\alpha \vee \chi_{k_1} \vee \chi_{k_1}^+ \vee \ldots \vee \chi_{k_1}^{+(M-k_1)}]$
\item[$(ii)$] $I_{k_j}$, $j = 2,\ldots,t$ denotes the interval
\[
[\alpha \vee \chi_{k_1} \vee \ldots \vee \chi_{k_j},\alpha \vee \chi_{k_1} \vee \chi_{k_2} 
        \vee \ldots \vee \chi_{k_j} \vee \chi_{k_{j-1}}^+ \vee \ldots \vee \chi_{k_{j-1}}^{(k_j - k_{j-1} - 1)}]
\]
\end{itemize}
\end{notation}
\begin{corollary}
\label{cor:multileveldecomposition}
Let $m < k_1 < k_2 < \ldots < k_t < M$. We have
\begin{eqnarray*}
|[\alpha,\beta]| & = &  
\sum_{\chi_{k_1}^{\ } \subseteq \PAB^{k_1}} |B_{k_1}| 
\sum_{\chi_{k_2} \subseteq \chi_{k_1}^{+(k_2-k_1)}} |I_{k-2}|
\sum_{\chi_{k_3} \subseteq \chi_{k_2}^{+(k_3-k_2)}} |I_{k-3}|
 \ldots \sum_{\chi_{k_t} \subseteq \chi_{k_{t-1}}^{+(k_t - k_{t-1})}} \ \  |I_{k_t}|\\
%|[\alpha,\beta]| & = & \sum_{\chi_{k_1}^{\ } \subseteq \PAB^{k_1}} \sum_{\chi_{k_2} \subseteq \chi_{k_1}^{+(k_2-k_1)}} \ldots \sum_{\chi_{k_t} \subseteq \chi_{k_{t-1}}^{+(k_t - k_{t-1})}} \\
%& & |[\alpha \vee \chi_{k_1},\alpha \vee \PAB^m \vee \ldots \PAB^{k_1-1} \vee \chi_{k_1}^{\ } ]| \\
%&  & .|[\alpha \vee \chi_{k_1} \vee \chi_{k_2},\alpha \vee \chi_{k_1} \vee \chi_{k_2} \vee \chi_{k_1}^+ \vee \ldots \chi_{k_1}^{+(k_2-k_1-1)}] \\
%& & \ldots  \\
%& & .|[\alpha \vee \chi_{k_1} \vee \ldots \chi_{k_p},\alpha \vee \chi_{k_1} \vee \ldots \chi_{k_p} \vee \chi_{k_{p-1}}^+ \vee \ldots \chi_{k_{p-1}}^{+(k_p-k_{(p-1)}-1)}]|\\
%& & \ldots  \\
%&  & .|[\alpha \vee \chi_{k_1} \vee \chi_{k_2} \vee \ldots \chi_{k_t},\alpha \vee \chi_{k_1} \vee \chi_{k_2} \vee \ldots \vee \chi_{k_t} \vee \chi_{k_t}^+ \vee \ldots \chi_{k_t}^{+(M-k_t)}]|
\end{eqnarray*}
\end{corollary}
\proof
By Theorem \ref{the:generalleveldecomposition} we have
\[
|[\alpha,\beta]| = \sum_{\chi_{k_1} \subseteq \PAB^{k_1}} |B_{k_1}|.|I_{k_1}|  
\]
Now apply Theorem \ref{the:generalleveldecomposition} to the interval $I_{k_1} = [\alpha \vee \chi_{k_1},\alpha \vee \chi_{k_1} \vee \chi_{k_1}^+ \vee \ldots \vee \chi_{k_1}^{+(M-k_1)}]$,
to get
\begin{eqnarray*}
|I_{k_1}| & =  \sum_{\chi_{k_2} \subseteq \PP_I^{k_2}} & |[\alpha \vee \chi_{k_1} \vee \chi_{k_2},\alpha \vee \chi_{k_1} \vee \chi_{k_1}^+ \ldots \vee \chi_{k_1}^{+(k_2-k_1-1)} \vee \chi_{k_2}]|.|I_{k_2}|
% & & .|[\alpha \vee \chi_{k_1} \vee \chi_{k_2}, \alpha \vee \chi_{k_1} \vee \chi_{k_2} \vee \chi_{k_2}^+ \vee \ldots \vee \chi_{k_2}^{+(M-k_2)}]|
\end{eqnarray*}
But  Lemma \ref{lem:lambdaischik} implies  $\PP_I^{k} = \chi_{k_1}^{+(k - k_1)}$ for any $k$, $k_1 < k \le k_2$, and hence
\begin{eqnarray*}
|I_{k_1}| & =  \sum_{\chi_{k_2} \subseteq \chi_{k_1}^{+(k_2-k_1)}} & |[\alpha \vee \chi_{k_1} \vee \chi_{k_2},\alpha \vee \chi_{k_1} \vee \chi_{k_1}^+ \ldots \vee \chi_{k_1}^{+(k_2-k_1-1)} \vee \chi_{k_2}]|.|I_{k_2}|
% & & .|[\alpha \vee \chi_{k_1} \vee \chi_{k_2}, \alpha \vee \chi_{k_1} \vee \chi_{k_2} \vee \chi_{k_2}^+ \vee \ldots \vee \chi_{k_2}^{+(M-k_2)}]|
\end{eqnarray*}
which proves the corollary for $t = 2$. Repeated application of Theorem  \ref{the:generalleveldecomposition}  on the intervals $I_{k_2},\ldots,I_{k_t}$ completes the proof.
\myendproof

Theorem \ref{the:generalleveldecomposition} is a generalization of Theorem \ref{the:evenoddleveldecomposition} 
that allows us to efficiently compute the size of an arbitrary interval by summing over selected levels not necessarily differing by two.
In contrast to Theorem  \ref{the:evenoddleveldecomposition}, the terms in the sum are not powers of two but products of sizes of intervals. For the size of these intervals,  no  closed form is presently known. However, the intervals of the form
\[
I_{k_p} = [\alpha \vee \chi_{k_1} \vee \ldots \chi_{k_p},\alpha \vee \chi_{k_1} \vee \ldots \chi_{k_{p}} \vee \chi_{k_{p-1}}^+ \ldots \vee \chi_{k_{p-1}}^{+(k_p-k_{p-1}-1)}]
\]
that show up as factors of terms in the expansion of Corollary \ref{cor:multileveldecomposition}, have an interesting structure 
that becomes apparent through the associated interval poset.
\begin{proposition}
\label{prop:horizontalintervals}
Let $m < k_1 <  k_2 \ldots < k_p < M$, and let 
\[
I \equiv I_{k_p} =  [\alpha \vee \chi_{k_1}^{\ } \vee \ldots \chi_{k_p}^{\ } ,\alpha \vee \chi_{k_1}^{\ }  \vee \ldots \chi_{k_p}^{\ }  \vee \chi_{k_{p-1}}^+ \ldots \vee \chi_{k_{p-1}}^{+(k_p-k_{p-1}-1)}]
\]
We have
\[
\PP_{I}^{\ } = \bigcup_{k_{p-1} < k < k_p} (\chi_{k_{p-1}}^{+(k-k_{p-1})} - \chi_{k_p}^{-(k_p - k)})
\]
and
\[
\PP_{I}^{k} = \chi_{k_{p-1}}^{+(k-k_{p-1})} - \chi_{k_p}^{-(k_p - k)}
\]
for $k_{p-1} < k < k_p$.
\end{proposition}
\proof
Since the antichain
$\chi_{k_{p-1}}^{+(k-k_{p-1})} - \chi_{k_p}^{-(k_p - k)}$ is $k$-uniform for each $k$, it suffices to show that
$\PP_{I}^k = \chi_{k_{p-1}}^{+(k-k_{p-1})} - \chi_{k_p}^{-(k_p - k)}$ for $k_{p-1} < k < k_p$, and empty otherwise.
Let $X \in \PP_{I}^{\ }$. 
From the definition of the interval poset, we infer that $\{X\} \not\le \alpha \vee \chi_{k_1}^{\ } \vee \ldots \chi_{k_p}^{\ }$ and $\{X\} \le \chi_{k_{p-1}}^+ \vee \ldots \chi_{k_(p-1)}^{+(k_p-k_{p-1}-1)}$.
Now let  $k_{p-1} < |X| = k <k_p$.
Notice that $\PP_{I}^{\ } \subseteq \PAB$ so that $|X| = k \Rightarrow X \in \PAB^k$.
By the definition of $()^+$, we now have $\{X\} \le  \chi_{k_{p-1}}^+ \vee \ldots \chi_{k_{p-1}}^{+(k_p-k_{p-1}-1)} \Rightarrow X \in  \chi_{k_{p-1}}^{+(k-k_{p-1})}$
and,  $\{X\} \not\le \chi_p \Rightarrow X \not\in \chi_p^{+(k_p - k)}$ and hence $X \in  \chi_{k_{p-1}}^{+(k-k_{p-1})} - \chi_{k_p}^{-(k_p - k)}$.
Conversely, we have $X \in  \chi_{k_{p-1}}^{+(k-k_{p-1})} \Rightarrow |X| > k_{p-1} \Rightarrow \{X\} \not\le \alpha \vee \chi_{k_1} \vee \ldots \chi_{k_{p-1}}$
and $X \not\in \chi_{k_p}^{-(k_p-k)} \Rightarrow \{X\} \not\le \chi_p$ such that $\{X\} \le  \chi_{k_{p-1}}^+ \vee \ldots \chi_{k_{p-1}}^{+(k_p-k_{p-1})} \Rightarrow X \in \PP_I^{\ }$ and we have proven
$k_{p-1} < k <k_p \Rightarrow \PP_I^k = \chi_{k_{p-1}}^{+(k-k_{p-1})} - \chi_{k_p}^{-(k_p - k)}$.

Now suppose $|X| \le k_{p-1}$ and $X \in \PP_I^{\ }$. Since $\{X\} \not\le \alpha \vee \chi_{k_1} \vee \ldots \chi_{k_{p}} \Rightarrow \{X\} \le  \chi_{k_{p-1}}^+ \vee \ldots \chi_{k_{p-1}}^{+(k_p-k_{p-1}-1)}$
and there is a set $Y \in \chi_{k_{p-1}}^{+(k-k_{p-1})}, X \subsetneq Y, k_{p-1}<k<k_p$.
Since $\{X\} \not\le \chi_{k_p}^{\ } \subseteq \chi_{k_{p-1}}^{+(k_p - k_{p-1})}$, we have $Y \in \chi_{k_{p-1}}^{+(k-k_{p-1})} -  \chi_{k_p}^{-(k_{p-1} -k)}$ and thus $Y \in \PP_I^{\ }$.
As a consequence of Theorem \ref{the:posetcharacteristic}, $X \subseteq Y \Rightarrow X \subsetneq X' \subsetneq Y, X' \in \PP_I^{\ }, |X'| = k_{p-1}$.
But any such set must satisfy $\{X'\} \le \chi_{k_{p-1}}^{+(k_p - k_{p-1})} \Rightarrow X' \in \chi_{k_{p-1}}^{\ } \Rightarrow \{X\} \le \chi_{k_{p-1}}^{\ }$ in contradiction with $X \in \PP_I^{\ }$.

If $|X| \ge k_p$  and $X \in \PP_I^{\ }$, then since the maximum size of a set in the top of the interval is $k_p$ (sets in $\chi_{k_p}^{\ }$) , we have $|X| = k_p$ and $\{X\} \le \chi_{k_p}^{\ }$ which is a contradiction.
This concludes the proof .
\myendproof
Note that, since $\chi_{k_p}^{\ } \subseteq \chi_{k_{p-1}}^{+(k_p - k_{p-1})}$, we always have $ \chi_{k_p}^{-(k_p - k)} \subseteq  \chi_{k_{p-1}}^{+(k-k_{p-1})}$.
\begin{remark}
An interesting case of Proposition \ref{prop:horizontalintervals} occurs when $k_p - k_{p-1} = 2$. In this case, the interval poset is given by
\[
\PP_I^{\ } = \chi_{k_{p-1}}^+ - \chi_{k_p}^- 
\]
and the sizes of all sets $X \in \PP_I^{\ }$ are equal to $k_{p-1} + 1$.
As a consequence, any subset of $\PP_I^{\ }$ is an antichain, and the number of antichains is given by 
\[
2^{|\PP_I^{\ }|} = 2^{|\chi_{k_{p-1}}^+| - | \chi_{k_p}^-|}
\]
This leads to an alternative proof of Theorem \ref{the:evenoddleveldecomposition}.
\end{remark}

% Corollary \ref{cor:multileveldecomposition} leads to Algorithm \ref{alg:multileveldecomposition}.

%\begin{algorithm}
%\caption{Given $k_1,\ldots,k_t$, $m < k_1 < k_2 < \ldots < k_t < M$, compute the size of the interval $[\alpha,\beta]$.
%\label{alg:multileveldecomposition}
%Use the notation $B_k$ and $I_k$ as in Corollary  \ref{cor:multileveldecomposition}.}
%\begin{algorithmic}
%\STATE{$Sum \gets 0$}
%\FOR {$\chi_{k_1} \in \PAB^{k_1}$}
%	\STATE{Compute $|B_{k_1}|$, through Algorithm \ref{alg:evenoddleveldecomposition}}
%	\STATE{$Sum \gets Sum + |B_{k_1}|*multilevel(\chi_{k_1}^{+(k_2 - k_1)},k_2,k_3,\ldots,k_t)$}
%\ENDFOR
%\RETURN{$Sum$}
%\end{algorithmic}
%\end{algorithm}

%\begin{algorithm}
%\caption{multilevel($\gamma$,$k_1,\ldots,k_t$)}
%\begin{algorithmic}
%\IF {$k_t == k_1$}
%\RETURN $|I_{k_t}|$
%\ENDIF
%\STATE{$Sum \gets 0$}
%\FOR {$\chi_{k_2} \in \gamma^{+(k_2 - k_1)}$}
%	\STATE{Compute $|I_{k_2}|$, through Algorithm \ref{alg:evenoddleveldecomposition}}
%	\STATE{$Sum \gets Sum + |I_{k_2}|*multilevel(\chi_{k_2}^{+(k_2 - k_1)},k_2,\ldots,k_t)$}
%\ENDFOR
%\RETURN{$Sum$}
%\end{algorithmic}
%\end{algorithm}

\section{Concluding remarks and future work}
We studied the structure of the set $\A_N$ of antichains of subsets of a finite set $N$.
Key elements in this study are the {\it intervals} with respect to a natural partial order on $\A_N$  based on
the operators {\it join} and {\it meet}. We introduced a new operator, the {\it direct product}, and
defined the {\it underlying poset} of intervals. 
The direct product is used in a Cartesian-like, decomposition theorem for general intervals (\ref{the:intervalintervaldecomposition}).
We derived a  characterisation for the underlying poset  (\ref{the:posetcharacteristic}), and used it in various expressions for the size of an interval based on unique decompositons of antichains as joins of {\it uniform} antichains (Section \ref{sec:decomposebyproperty}).
%\ref{the:evenoddleveldecomposition}, \ref{the:generalleveldecomposition}, \ref{cor:multileveldecomposition}, \ref{prop:horizontalintervals}).
The concept of a {\it largest non-dominating antichain} is used in Theorem  (\ref{the:largestnondominatingdecomposition}), which establishes a partition of a general interval associated with a general antichain.
Finally, we derived a {\it direct-join decomposition} of intervals with a specific structure.

Overall, we have presented a systematic analysis of the interval structure of $\A_N$ opening perspectives for coming work.
In the past \cite{WIEDEMANN, BAKOEV, FYDITEK} a specific sum formula for the Dedekind numbers, $|\A_N|$ in our notation, has been used to compute the largest cases.
The sum formula relies on the sizes of specific intervals.
A description of such sum formulas can be found in \cite{PDCSDW2013}. 
The properties developed in the present paper allow to derive a range of such sum formulas which we will study in coming work.

On the structural side, the underlying poset and the direct join decomposition provide a new perspective for the study of intervals.
In particular, isomorphism between intervals and reduction of intervals to a canonical form, can be tackled using these constructs.
Such a canonical form would automatically lead to a canonical form for the antichains which are the borders of the intervals.

Properties of antichains such as connectedness, the fact that two elements always come together and symmetry under permutation allow recursion formulas to be derived. 
So far, we did not come up with  families of antichains that can both be used in a recursion formula and are asymptotically smaller than $|\A_N|$.
The properties of intervals are a tool to study such families, and we plan to continue in this direction.

\section{Acknowledgement}
Work supported by the Belgian Science Policy Office (BELSPO) in the Interuniversity Attraction Pole COMEX. (http://comex.ulb.ac.be)

\section*{References}
\bibliography{mybibfile}

\begin{thebibliography}{15}
\providecommand{\natexlab}[1]{#1}
\providecommand{\url}[1]{\texttt{#1}}
\expandafter\ifx\csname urlstyle\endcsname\relax
  \providecommand{\doi}[1]{doi: #1}\else
  \providecommand{\doi}{doi: \begingroup \urlstyle{rm}\Url}\fi

\bibitem[Bakoev(2012)]{BAKOEV}
Valentin Bakoev.
\newblock One more way for counting monotone boolean functions.
\newblock In \emph{Thirteenth International Workshop on Algebraic and
  Combinatorial Coding Theory}, page 47Ð52, 2012.

\bibitem[Church(1940)]{CHURCH1940}
R.~Church.
\newblock Numerical analysis of certain free distributive structures.
\newblock \emph{Duke. Math. J.}, 6:\penalty0 732Ð734, 1940.

\bibitem[Church(1965)]{CHURCH1965}
R~Church.
\newblock Enumeration by rank of the elements of the free distributive lattice
  with 7 generators.
\newblock \emph{Notices Amer. Math. Soc.}, 12:\penalty0 724, 1965.

\bibitem[De~Causmaecker and De~Wannemacker(2011)]{PDCSDW2011}
Patrick De~Causmaecker and Stefan De~Wannemacker.
\newblock Partitioning in the space of anti monotonic functions.
\newblock \emph{arXiv:1103.2877 [math.NT]}, 2011.

\bibitem[De~Causmaecker and De~Wannemacker(2013)]{PDCSDW2013}
Patrick De~Causmaecker and Stefan De~Wannemacker.
\newblock Decomposition of intervals in the space of anti-monotonic functions.
\newblock In M.~Ojeda-Aciego and J.~Outrata, editors, \emph{Proceedings of the
  Tenth International Conference on Concept Lattices and Their Applications
  (CLA 2013)}, volume 1062, pages 57--67, La Rochelle, France, 15-18 October
  2013 2013. CLA 2013.

\bibitem[Dedekind(1897)]{DEDEKIND}
Richard Dedekind.
\newblock Uber zerlegungen von zahlen durch ihre grossten gemeinsamen teiler.
\newblock \emph{Gesammelte Werke 2}, pages 103--148, 1897.

\bibitem[Fidytek et~al.(2001)Fidytek, Mostowski, Somla, and
  Szepietowski]{FYDITEK}
Robert Fidytek, Andrzej~W. Mostowski, Rafal Somla, and Andrzej Szepietowski.
\newblock Algorithms counting monotone boolean functions.
\newblock \emph{Information Processing Letters}, 79:\penalty0 203Ð209, 2001.

\bibitem[Kahn(2002)]{KAHN}
Jeff Kahn.
\newblock Entropy, independent sets and antichains: a new approach to
  dedekind's problem.
\newblock \emph{Proc. Amer. Math. Soc.}, 130 (2):\penalty0 371--378, 2002.

\bibitem[Kleitman(1969)]{KLEITMAN}
Daniel Kleitman.
\newblock On dedekind's problem: the number of isotone boolean functions.
\newblock \emph{Proc. Amer. Math. Soc.}, 21:\penalty0 677--682, 1969.

\bibitem[Kleitman and Markowsky(1975)]{KLEITMAN_MARKOWSKY}
Daniel Kleitman and George Markowsky.
\newblock On dedekind's problem: the number of isotone boolean functions. ii.
\newblock \emph{Trans. Amer. Math. Soc.}, 213:\penalty0 373--390, 1975.

\bibitem[Korshunov(1981)]{KORSHUNOV}
Aleksej~Dmitrievich Korshunov.
\newblock The number of monotone boolean functions (russian).
\newblock \emph{Problemy Kibernet.}, 38:\penalty0 5--108, 1981.

\bibitem[Sloane()]{SLOANE}
N.J.A. Sloane.
\newblock The on-line encyclopedia of integer sequences. (oeis).
\newblock \emph{http://www.research.att.com/~njas/sequences/}.

\bibitem[Stephen and Yusun(2014)]{TAMON2014}
Tamon Stephen and Timothy Yusun.
\newblock Counting inequivalent monotone boolean functions.
\newblock \emph{Discrete Appl. Math.}, 167:\penalty0 15--24, April 2014.
\newblock ISSN 0166-218X.
\newblock \doi{10.1016/j.dam.2013.11.015}.
\newblock URL \url{http://dx.doi.org/10.1016/j.dam.2013.11.015}.

\bibitem[Ward(1946)]{WARD}
M.~Ward.
\newblock Note on the order of free distributive lattices.
\newblock \emph{Bull. Amer. Math. Soc.}, 52:\penalty0 423, 1946.

\bibitem[Wiedemann(1991)]{WIEDEMANN}
Doug Wiedemann.
\newblock A computation of the eighth dedekind number.
\newblock \emph{Order}, 8 (1):\penalty0 5--6, 1991.

\end{thebibliography}

\end{document}